\documentclass[a4paper,12pt,reqno]{article}
\usepackage[T1]{fontenc}
\usepackage[hmargin=1in,vmargin=1in,nohead]{geometry}
\usepackage{textcomp}
\usepackage{amsmath}
\usepackage{graphicx,epsfig,color}
\usepackage{amssymb}
\usepackage{esint}

\makeatletter

\newtheorem{thm}{Theorem}
\newtheorem{lem}[thm]{Lemma}

\newtheorem{prop}[thm]{Proposition}
\newtheorem{rem}{Remark}[section]

\newtheorem{Example}{Example}[section]
\newtheorem{Def}{Definition}
\DeclareMathOperator{\cg}{\textbf{[}}
\DeclareMathOperator{\cd}{\textbf{]}}
\newcommand{\cro}[1]{\cg {#1} \cd}
\date{\today}
\title{ \bf Distinguished exchangeable coalescents and generalized Fleming-Viot processes with immigration}
\author{ \bf Cl\'ement Foucart \\
\\
\emph {Laboratoire de Probabilit\'es et Mod\`eles Al\'eatoires} \\
\emph {Universit\'e Pierre et Marie Curie}\\ 
\emph{4 Place Jussieu- 75252 Paris Cedex 05- France}}
\begin{document}
\maketitle{}
\begin{center}
DRAFT VERSION
\end{center}

\begin{abstract} 
Coalescents with multiple collisions (also called $\Lambda$-coalescents or simple exchangeable coalescents) are used as models of genealogies. We study a new class of Markovian coalescent processes connected to a population model with immigration. Imagine an infinite population with immigration labelled at each generation by $\mathbb{N}:=\{1,2,...\}$. Some ancestral lineages cannot be followed backwards after some time because their ancestor is outside the population. The individuals with an immigrant ancestor constitute a distinguished family and we define exchangeable distinguished coalescent processes as a model for genealogy with immigration, focussing on simple distinguished coalescents, i.e such that when a coagulation occurs all the blocks involved merge as a single block. These processes are characterized by two finite measures on $[0,1]$ denoted by $M=(\Lambda_{0},\Lambda_{1})$. We call them $M$-coalescents. We show by martingale arguments that the condition of coming down from infinity for the $M$-coalescent coincides with that obtained by Schweinsberg for the $\Lambda$-coalescent. In the same vein as Bertoin and Le Gall, $M$-coalescents are associated with some stochastic flows. The superprocess embedded can be viewed as a generalized Fleming-Viot process with immigration. The measures $\Lambda_{0}$ and $\Lambda_{1}$ specify respectively the reproduction and the immigration. The coming down from infinity of the $M$-coalescent will be interpreted as the initial types extinction: after a certain time all individuals are immigrant children.
\end{abstract}
 \vspace{9pt} \noindent {\bf Key words.}
{Exchangeable partitions}, {coalescent theory}, {genealogy for a population with immigration},
{stochastic flows}, {coming down from infinity.}
\par \vspace{9pt} \noindent {\bf e-mail.} {clement.foucart@etu.upmc.fr}
\par

\newpage
\newpage
\section{Introduction}\label{secIntro}
Pitman, \cite{Pitman2} and Sagitov, \cite{Sagitov} defined in 1999 the class of $\Lambda$-coalescent processes, sometimes also called simple exchangeable coalescents. These coalescent processes appear as models for the genealogy of certain haploid populations with fixed size. The general motivation of this work is to define a new class of coalescent processes that may be used to describe the genealogy of a population with immigration. Heuristically, let us imagine an infinite haploid population with immigration described at each generation by $\mathbb{N}:=\{1,2,3,...\}$. This means that each individual has at most one parent in the population at the previous generation; indeed, immigration implies that some individuals may have parents outside this population (they are children of immigrants). Sampling $n$ individuals in the population at some fixed generation, we group together the individuals with the same parent at the preceding generation. The individuals with an immigrant parent constitute a special family. We get a partition of $\mathbb{N}$ where each block is a family. \\
To give to the population a full genealogy, we may imagine a generic external ancestor, say $0$, to distinguish the immigrants family. This way, all families will have an ancestor at the preceding generation. Following the ancestral lineage of an individual backwards in time, it may coalesce with some others in $\mathbb{N}$ or reach $0$. In that last case the lineage is absorbed at $0$. We call $0$ the immigrant ancestor, and we shall therefore work with partitions of $\mathbb{N}\cup \{0\}=\mathbb{Z}_{+}$. We view the block containing $0$ as distinguished and then we speak of distinguished partitions. As usual, a partition is identified with the sequence of its blocks in the increasing order of their smallest element. The distinguished block is thus the first.\\
For a population with no immigration, Kingman introduced exchangeable random partitions of $\mathbb{N}$. A random partition is exchangeable if and only if its law is invariant under the action of permutations of $\mathbb{N}$. The distinguished partitions appearing in our setting are not exchangeable on $\mathbb{Z}_{+}$, however their laws are invariant under the action of permutations $\sigma$ of $\mathbb{Z}_{+}$ such that $\sigma(0)=0$. These partitions are called exchangeable distinguished partitions. We will present an extension of Kingman's theorem that determines their structure via a paint-box construction.
\\
This allows us to define, following the approach in Bertoin's book \cite{coursbertoin}, a new class of coalescent processes, which we call exchangeable distinguished coalescents. An exchangeable distinguished coalescent is characterized in law by a measure $\mu^{0}$ on the space of partitions of $\mathbb{Z}_{+}$, called the distinguished coagulation measure. The extension of Kingman's theorem enables us to characterize this measure, and when $\mu^{0}$ is carried on the subset of simple distinguished partitions (which have only one non-trivial block), we get a representation involving two finite measures on $[0,1]$: $M=(\Lambda_{0},\Lambda_{1})$. We call $M$-coalescents this sub-class of distinguished coalescents. The restriction of an $M$-coalescent to each finite subset containing $0$, is a Markovian coalescent chain with the following transition rates: when the partition restricted to $\mathbb{N}$ has $b$ blocks, two kinds of jumps are allowed : for $b\geq k\geq 2$ each $k$-tuple of blocks not containing $0$ can merge to form a single block at rate $\int_{0}^{1}x^{k-2}(1-x)^{b-k}\Lambda_{1}(dx)$, and for $b\geq k\geq 1$, each $k$-tuple of blocks not containing $0$ can merge with the one containing $0$ at rate $\int_{0}^{1}x^{k-1}(1-x)^{b-k}\Lambda_{0}(dx)$.
\\ 
Next, we study a classical question for coalescent processes : a coalescent process starting from infinitely many blocks, is said to come down from infinity if its number of blocks instantaneously becomes finite. An interesting result is that the condition for $M$-coalescents to come down does not depend on $\Lambda_{0}$ and is the same for the $\Lambda_{1}$-coalescent found by Schweinsberg \cite{CDI}. 
\\
In the last section, we define some stochastic flows connected with $M$-coalescents. The model of continuous population embedded in the flow can be viewed as a generalized Fleming-Viot process with immigration. As in \cite{LGB1}, the stochastic flows involved allow us to define simultaneously a population model forward in time and its genealogical process backward in time. A duality between $M$-generalized Fleming-Viot processes with immigration and $M$-coalescents will be studied.
\\
\\
In a forthcoming paper, we will give a different approach to construct the generalized Fleming-Viot processes with immigration by introducing some stochastic flows of partitions. Our method will draw both on the works of Donnelly-Kurtz \cite{DonnKurtz} and of Bertoin-Le Gall \cite{LGB1}. Some ideas of Birkner \textit{et al} in \cite{Birk} may be applied to establish a link between certain branching processes with immigration and $M$-generalized Fleming-Viot processes with immigration.
\\
\\
\textbf{Outline.} The paper is organized as follows. In Section 2, we recall some basic facts on random partitions, and we give some fundamental properties of exchangeable distinguished partitions (about existence of asymptotic frequencies, paint-box representation). In Section 3, we define exchangeable distinguished coalescents. We establish a characterization of their laws by an exchangeable measure $\mu^{0}$ on the space of the distinguished partitions. The structure of $\mu^{0}$ is entirely described which enables us to study the dust. The main reference is Chapters 2 and 4 of Bertoin's book \cite{coursbertoin}. The construction of exchangeable distinguished coalescents is very close to that for exchangeable coalescents of \cite{coursbertoin}. In Section 4, we focus on $M$-coalescents and study the coming down from infinity. In particular, our approach provides a new proof of Schweinsberg's result, see \cite{CDI}, about necessary and sufficient conditions to come down from infinity for $\Lambda$-coalescents based on martingale arguments. In Section 5, we introduce certain stochastic flows encoding $M$-coalescents. As in \cite{coursbertoin} and \cite{LGB1}, these flows allow us to define a population model with immigration called $M$-generalized Fleming-Viot process with immigration.   
\section{Distinguished partitions}
We begin with some general notation and properties which we will use constantly in the following sections.\\
For every integer $n\geq 1$, we denote by $[n]$ the set $\{1,...,n\}$ and by $\mathcal{P}_{n}$ the set of its partitions. The set of partitions of $\mathbb{N}$ is denoted by $\mathcal{P}_{\infty}$. Let $\pi \in \mathcal{P}_{\infty}$, we identify the set $\pi$ with the sequence $(\pi_{1},\pi_{2},...)$ of the blocks of $\pi$ enumerated in increasing order of their least element: for every $i\leq j$, $\min \pi_{i} \leq \min \pi_{j}$. The number of blocks of $\pi$ is $\# \pi$. 
For all $\pi \in \mathcal{P}_{\infty}$ and $n\in \mathbb{N}$, $\pi_{|[n]}\in \mathcal{P}_{n}$ is by definition the restriction of $\pi$ to $[n]$. We denote by $\mathcal{P}_{\textbf{m}}$ the set of mass-partitions, meaning the decreasing sequences with sum less than or equal to $1$: 
$$\mathcal{P}_{\textbf{m}}:=\{s=(s_{1},s_{2},...); \sum_{i\geq 1} s_{i} \leq 1, s_{1}\geq s_{2}\geq ...\geq 0\}.$$ 
Given a partition $\pi=(B_{1},B_{2},...)$ and a block $B$ of that partition, we say that $B$ has an asymptotic frequency, denoted by $|B|$, if the following limit exists: $$|B|:=\underset{n\rightarrow \infty}{\lim} \frac{\#(B\cap[n])}{n}.$$ If each block of a partition has asymptotic frequency, this partition is said to have asymptotic frequencies. For $\pi \in \mathcal{P}_{\infty}$ possessing asymptotic frequencies, $|\pi|^{\downarrow}$ is the mass partition associated with $\pi$ that is $(|\pi|^{\downarrow}_{i})_{i\in \mathbb{N}}$ is the rearrangement in decreasing order of $(|\pi_{i}|)_{i\in \mathbb{N}}$. 
For every $n\in \mathbb{N}$, a permutation of $[n]$ is a bijection $\sigma:[n]\mapsto [n]$. For $n=\infty$, we define a permutation of $\mathbb{N}$ to be a bijection $\sigma$ of $\mathbb{N}$ such that $\sigma(k)=k$ when $k$ is large enough. We define the equivalence relation $\underset{\pi}{\sim}$ by $i\underset{\pi}{\sim} j$ if $i$ and $j$ are in the same block of $\pi$. We denote $\sigma \pi$ the partition defined by
\begin{center}
$i\underset{\sigma \pi}{\sim} j \Longleftrightarrow \sigma(i)\underset{\pi}{\sim} \sigma(j).$
\end{center}
We stress that due to the ranking of the blocks, $(\sigma \pi)_{i} = \sigma^{-1}(\pi_{\eta(i)})$ for a certain permutation $\eta$.
\\
A random partition $\pi$ of $\mathbb{N}$ is exchangeable if $\sigma \pi$ and $\pi$ have the same law, for every permutation $\sigma$ of $\mathbb{N}$. Kingman established a correspondence between exchangeable partitions laws and mass-partitions via the paint-box partitions. We recall briefly the construction of paint-boxes. Let $s$ be an element of $\mathcal{P}_{\textbf{m}}$. Let $\mathcal{V}$ be an open subset of $]0,1[$ such that the ranked sequence of lengths of its interval components is given by $s$. Let $U_{1},...$ be an i.i.d sequence of uniform variables on $[0,1]$. A $s$-paint-box is the partition $\pi$ induced by the equivalence relation: 
\begin{center} $\forall i\neq j$: $i \underset{\pi}{\sim} j \Longleftrightarrow U_{i}$ and $U_{j}$ belong to the same interval component of $\mathcal{V}$.\end{center} Kingman proved that any exchangeable partition is a mixture of paint-boxes. We denote by $\rho_{s}$ the law of a $s$-paint-box. 
\\
\\
As explained in the Introduction, we now extend this setting by distinguishing a block, working with partitions of $\mathbb{Z}_{+}$.
\begin{Def} 
A distinguished partition $\pi$ is a partition of $\mathbb{Z}_{+}$ where the block containing $0$ is viewed as a distinguished block. Ranking the blocks in the order of their least element, the first block $\pi_{0}$ contains $0$ and is the distinguished block of $\pi$.
\end{Def}
We denote by $\cro{n}$ the set $\{0,1,...,n\}$, $\mathcal{P}_{n}^{0}$ is the space of distinguished partitions of $\{0,...,n\}$. For $n=\infty$, we agree that $\cro{\infty}=\mathbb{Z}_{+}=\{0,1,...\}$ and then $\mathcal{P}_{\infty}^{0}$ is the space of partitions of $\mathbb{Z}_{+}$. A first basic property is the compactness of the space $\mathcal{P}^{0}_{\infty}$ for the distance defined by $$d(\pi,\pi')=(1+\max\{n\geq 0, \pi_{|\cro{n}}=\pi'_{|\cro{n}}\})^{-1}.$$ See \cite{coursbertoin} for a proof. Let $\pi \in \mathcal{P}^{0}_{n}$, for all $n'\in \cro{\infty}$ such that $n'\geq n$, we define $$\mathcal{P}^{0}_{n',\pi}=\{\pi'\in \mathcal{P}^{0}_{n'}; \pi'_{|\cro{n}}=\pi\}.$$
\\
A random distinguished partition is a random element of $\mathcal{P}^{0}_{\infty}$ equipped with the $\sigma$-field generated by the finite unions of the sets $\mathcal{P}^{0}_{n,\pi}$ (that corresponds to the Borelian $\sigma$-field for $d$).
\\
\\
In the same way, we introduce the set of distinguished mass-partitions, meaning the sequences $s=(s_{i})_{i\geq 0}$ of non-negative real numbers such that $\sum_{i\geq 0}s_{i}\leq 1$, ranked in decreasing order apart from $s_{0}$ 
$$\mathcal{P}^{0}_{\textbf{m}}:=\{s=(s_{0},s_{1},...); \sum_{i\geq 0} s_{i} \leq 1,  s_{0}\geq 0, s_{1}\geq s_{2}\geq ...\geq 0\}.$$ 
We identify the sets $\mathcal{P}_{\textbf{m}}$ and $\{s\in \mathcal{P}^{0}_{\textbf{m}}; s_{0}=0\}$. The dust of $s$ is by definition the quantity $\delta:=1-\sum_{i=0}^{\infty}s_{i}$. A (distinguished) mass-partition is said to be improper if the dust is positive. For $\pi \in \mathcal{P}^{0}_{\infty}$ having asymptotic frequencies, $|\pi|^{\downarrow}$ is the distinguished mass-partition associated with $\pi$ that is $|\pi |^{\downarrow}_{0}=|\pi_{0}|$ and $(|\pi|^{\downarrow}_{i})_{i\in \mathbb{N}}$ is the rearrangement in decreasing order of $(|\pi_{i}|)_{i\in \mathbb{N}}$. We stress that by definition $|\pi|^{\downarrow}_{0}=|\pi_{0}|$.\\
We define a permutation of $\mathbb{Z}_{+}$ to be a bijection $\sigma$ of $\mathbb{Z}_{+}$ such that $\sigma(k)=k$ when $k$ is large enough. Note that any permutation of $\mathbb{N}$ can be extended to a permutation of $\mathbb{Z}_{+}$ by deciding that $\sigma(0)=0$.
\begin{Def} 
A random distinguished partition $\pi$ is exchangeable if $\sigma \pi$ and $\pi$ have the same law for every permutation $\sigma$ of $\mathbb{Z}_{+}$ such that $\sigma(0)=0$.
\end{Def}
It is easily seen that the restriction of an exchangeable distinguished partition to $\mathbb{N}$ is exchangeable. The converse may fail: there exist distinguished partitions which are not exchangeable though their restriction to $\mathbb{N}$ is exchangeable. We construct a counter-example: let $\pi$ be a non-degenerate exchangeable random partition and $\pi^{0}$ obtained from $\pi$ by distinguishing the block containing $1$, i.e $\pi^{0}=(\pi_{1}\cup \{0\},\pi_{2},...)$ with blocks enumerated in order of appearance. The restriction $\pi^{0}_{|\mathbb{N}}=\pi$ is exchangeable. The structure of $\pi$ implies that $\mathbb{P}[\pi^{0}_{|\cro{2}}=(\{0,1\},\{2\})]=\mathbb{P}[1 \not \underset{\pi}{\sim} 2]>0$. Let $\sigma$ be the permutation of $\cro{2}$: $\sigma (0)=0, \sigma (1)=2, \sigma (2)=1$. We have $\mathbb{P}[\pi^{0}_{|\cro{2}}=(\{0,\sigma(1)\},\{\sigma(2)\})]=\mathbb{P}[\pi_{1}=\{2\}]=0$. We thus found a permutation such that $\mathbb{P}[\pi^{0}=\sigma(\pi^{0}_{0},...)]\neq  \mathbb{P}[\pi^{0}=(\pi^{0}_{0},...)]$. 
\\
\\
We define now the distinguished paint-boxes and extend the Kingman's correspondence to exchangeable distinguished partitions.
\begin{Def} \label{paintbox}
A distinguished paint-box can be constructed in the following way:\\
Let $s$ be a distinguished mass partition, we denote by $\delta$ its dust. Denote by $\partial$, an element which does not belong to $\mathbb{Z}_{+}$. Let $\xi$ be a probability on $\mathbb{Z}_{+} \cup \{ \partial \}$ such that for all $k\geq 0, \xi(k)=s_{k}$ and $\xi(\partial)=\delta$. Drawing $X_{1},X_{2},...$ a sequence of i.i.d random variables with distribution $\xi$ and $X_{0}=0$. 
A $s$-distinguished paint box is defined by : $\forall i\neq j \geq 0$,
\begin{center}
$\Pi^{0}:$ $ i \sim j \Longleftrightarrow X_{i}=X_{j}\neq \partial$.
\end{center}
In particular, $\Pi^{0}_{0}:=\{i\geq 0; X_{i}=0\}=\{i\geq 1; X_{i}=0\}\cup \{0\}$.
\end{Def}
We denote by $\rho^{0}_{s}$ the law of an $s$-distinguished paint box. When $s_{0}=0$, the block $\Pi^{0}_{0}$ is the singleton $\{0\}$ and the $s$-distinguished paint box restricted to $\mathbb{N}$ is a classical $s$-paint box partition of $\mathbb{N}$. According to the previous notation for paint-boxes, we will denote its law by $\rho_{s}$. 
\\
\\
Another way to define an $s$-distinguished paint-box is to consider a sub-probability $\alpha$ on $[0,1[$ and set $s=(\alpha(0), \alpha(x_{1}), \alpha(x_{2}),...)$ where $x_{1},x_{2},...$ are the atoms of $\alpha$ in $]0,1[$ ranked in decreasing order of their sizes. Let $X_{1},X_{2}...$ be independent with law $\alpha$ and $X_{0}=0$, the partition $\Pi^{0}$ defined by: $i\neq j$: $i \sim j \Longleftrightarrow X_{i}=X_{j}$ is an $s$-distinguished paint-box. 
\\
Equivalently, we can work with uniform variables: an interval representation of $s$ is a collection of disjoint intervals $(A_{0},A_{1},A_{2},...)$, where $A_{0}$ is $[0,s_{0}]$, and $(A_{i})_{i\geq 1}$ such that the decreasing sequence of their lengths is $(s_{1},s_{2},...)$. If we draw an infinite sequence of uniform independent variables $(U_{i})_{i\geq 1}$ and fix $U_{0}=0$. The partition of $\mathbb{Z}_{+}$ defined by: $\pi^{0}:=i\sim j$ if and only if $U_{i}$ and $U_{j}$ fall in the same interval (if $U_{i}$ falls in the dust of $s$ then $\{i\}$ is a singleton block of $\pi^{0}$) is a $s$-distinguished paint box. We stress that its law does not depend on the choice of intervals $(A_{1},A_{2}...)$ and then we can choose $A_{i}:=[s_{0}+...+s_{i-1},s_{0}+...+s_{i}[$ for all $i\geq 1$.
\begin{prop} Let $s\in \mathcal{P}^{0}_{\textbf{m}}$ and $\pi^{0}$ a $s$-distinguished paint-box. 
\begin{itemize}
\item[i)]  The distinguished paint-box $\pi^{0}$ is exchangeable.
\item[ii)] $\pi^{0}$ has asymptotic frequencies, and more precisely $|\pi^{0}|^{\downarrow}=s$.
\item[iii)]For every $i\in \mathbb{Z}_{+}$, if $|\pi^{0}_{i}|=0$ then $\pi^{0}_{i}$ is a singleton or empty.
\item[iv)] $s$ is improper if and only if some blocks different from $\pi^{0}_{0}$ are singletons. In that case the set of singletons $\{i\in \mathbb{Z}_{+}: i$ is a singleton of $\pi^{0}\}$ has an asymptotic frequency given by the dust $\delta=1-\sum_{i=0}^{\infty}s_{i}$ a.s
\item[v)]  We have $\rho^{0}_{s}(0$ is singleton)$=0$ if $s_{0}>0$ and $1$ otherwise, and $\rho^{0}_{s}(1, 2,...,q$ are singletons)=$\delta^{q}$, for $q\geq 1$.\end{itemize}
\end{prop}
\textit{Proof.} See the proof of Proposition 2.8 in \cite{coursbertoin}. $\square$
\\
It remains to see whether the distinguished paint-box construction of Definition \ref{paintbox} yields all the exchangeable distinguished partitions. 
\begin{thm} \label{repres} Let $\pi^{0}$ be a random distinguished partition. The following assertions are equivalent:
\begin{itemize}
\item[i)]  $\pi^{0}$ is exchangeable
\item[ii)] There exists a random distinguished mass-partition $S=(S_{0},S_{1},...)$ such that conditionally given $S=s$, $\pi^{0}$ has the law of an $s$-distinguished paint box $(\rho^{0}_{s})$. Further $|\pi^{0}|^{\downarrow}=S$.
\end{itemize}
\end{thm}
\textit{Proof.} A mixture of distinguished paint-boxes is still exchangeable, this shows that \textit{ii)} implies \textit{i)}. Let $\pi^{0}$ be exchangeable. We adapt a proof of Aldous \cite{Aldous}, see also \cite{coursbertoin} p101. We call selection map, any random function $b:\mathbb{Z}_{+} \rightarrow \mathbb{Z}_{+}$ that maps all the points of the block $\pi^{0}_{0}$ to $0$, and all the points of a block $\pi^{0}_{i}$ for $i\geq 1$ to the same point of that block. Let $U_{0}=0$, $(U_{i})_{i\geq 1}$ be i.i.d uniform on $[0,1]$, independent of $\pi^{0}$ and of the selection map $b$. We define $X_{n}=U_{b(n)}$. The law of $(X_{n}, n\geq 1)$ does not depend on the choice of $b$. The key of the proof is the exchangeability of $(X_{n})_{n\geq 1}$. Let $\sigma$ be a permutation with $\sigma(0)=0$, we have $$X_{\sigma(n)}=U_{b(\sigma(n))}=U'_{b'(n)},$$
where $U'_{i}=U_{\sigma(i)}$ and $b'=\sigma^{-1}\circ b \circ \sigma$. We verify that $b'$ is a selection map for the partition $\sigma \pi^{0}$. Let $i\geq 0$, and $n\in \sigma \pi^{0}_{i}$, by definition of $\sigma \pi^{0}$, $\sigma\pi^{0}_{i}=\sigma^{-1}(\pi^{0}_{\eta(i)})$ for a certain permutation $\eta$ such that $\eta(0)=0$, then there exists $k\in \pi^{0}_{\eta(i)}$ such that $n=\sigma^{-1}(k)$. For $i=0$, $k\in \pi^{0}_{0}$ and then $b'(n)=\sigma^{-1}(b(k))=\sigma^{-1}(0)=0$ for all $n\in \sigma \pi^{0}_{0}$. For $i\geq 1$, we clearly have that $b'(n)=\sigma^{-1}\circ b(k)$ depends only on $i$.
\\
The sequence $(U'_{i}, i\geq 1)$ has the same law as $(U_{i}, i\geq 1)$. By exchangeability and independence of $(U_{j})$ and $\pi^{0}$: $((U'_{n})_{n\geq 1}, \sigma \pi^{0})$ has the same law as $((U_{n})_{n\geq 1}, \pi^{0})$, and the sequence $(X_{n},n \geq 1)$ is exchangeable. By the de Finetti theorem, we have that conditionally on the random probability measure $\mu:=\underset{n\rightarrow \infty}{\lim} \frac{1}{n} \sum_{i=1}^{n} \delta_{X_{i}}$, $(X_{n}, n\geq 1)$ are i.i.d random variables with distribution $\mu$. Moreover, by the definition of $X_{n}$, $i\underset{\pi}{\sim} j$ if and only if $X_{i}=X_{j}$. We deduce that conditionally given $\mu=\alpha$, the distinguished partition $\pi^{0}$ is a $s(\alpha)$-distinguished paint box with $s(\alpha):=(\alpha(0),\alpha(x_{1}),...)$. By the distinguished paint-box construction, on $\{\mu=\alpha\}$, $s(\alpha)$ is the mass-partition of $\pi^{0}$ and so $|\pi^{0}|^{\downarrow}=s(\mu)$, moreover the random sequence $S:=s(\mu)$ in $\mathcal{P}^{0}_{\textbf{m}}$ verifies the assertion \textit{ii)}. To conclude the random partition $\pi^{0}$ has the law of a $\eta$-mixture of distinguished paint-boxes, where $\eta$ is the law of $|\pi^{0}|^{\downarrow}$.
$\square$
\\
\\
This theorem sets up a bijection between probability distributions for exchangeable distinguished partitions and probability distributions on the space of distinguished mass-partitions, $\mathcal{P}^{0}_{\textbf{m}}$: $$\mathbb{P}[\pi^{0}\in .]=\int_{\mathcal{P}^{0}_{m}}\rho^{0}_{s}(.)\mathbb{P}(|\pi^{0}|^{\downarrow}\in ds).$$
\begin{rem}
Let $\pi^{0}$ be an exchangeable distinguished partition, as for exchangeable partitions, see \cite{courspitman}, we can show that there exists a function $p$ such that $$\mathbb{P}[\pi^{0}_{|\cro{n}}=(B_{0},...,B_{k})]=p(n_{0},...,n_{k}).$$ where $n_{i}=\#B_{i}$, $i=0,...,k$. Contrary to exchangeable random partitions, the function $p$ is not totally symmetric but only invariant by permutations of the arguments $(n_{1},...,n_{k})$. Indeed by exchangeability $\mathbb{P}[\pi^{0}=(\sigma B_{0},...,\sigma B_{k})]=p(\#B_{0},...,\#B_{k})=p(\#B_{\eta(0)},...,\#B_{\eta(k)})$ where $\eta$ is the permutation such that $\sigma \pi^{0}_{i}=\sigma^{-1}(\pi^{0}_{\eta(i)})$. Due to the assumption $\sigma(0)=0$, the permutation $\eta$ is such that $\eta(0)=0$. The exchangeable distinguished partitions are thus special cases of partially exchangeable partitions in the sense of \cite{Pitman1}. \\
We mention that Donnelly and Joyce, \cite{Donnelly} defined the exchangeable ordered partitions for which all the blocks are distinguished (they speak about exchangeable random ranking). They obtain a Kingman's representation for the exchangeable ordered partition structure. For every $\mu$ probability on $[0,1]$, an "ordered paint-box" is constructed from a sequence of i.i.d. $\mu$ variables. We stress that, contrary to an exchangeable distinguished paint-box, the law of a $\mu$-ordered paint-box depends on \textbf{the order of the atoms} of $\mu$. Exchangeable ordered partitions are also partially exchangeable but the function $p$ has no symmetry properties. Gnedin in \cite{Gnedin} gives a representation of exchangeable compositions which are a generalization of exchangeable random rankings.\\
We could define distinguished partitions with several distinguished blocks. It corresponds to a population with several sources of immigration: each distinguished block gathers the progeny of an immigration source. For sake of simplicity, we distinguish here just one block.
\end{rem}
In the next section, we define distinguished coalescents which can be interpreted as a genealogy for a population with immigration. 
\section{Distinguished coalescents}
Imagine an infinite haploid population with immigration. We denote by $\Pi^{0}(t)$ the partition of the current population into families having the same ancestor $t$ generations earlier. As explained in the Introduction, individuals who have no ancestor in the population at generation $t$, form the distinguished block of $\Pi^{0}(t)$. Actually, the individual $0$ can be viewed as their common ancestor. When some individuals have the same ancestor at a generation $t$, they have the same ancestor at any generation $t'\geq t$. In terms of partitions, all integers in the same block of $\Pi^{0}(t)$, are in the same block of $\Pi^{0}(t')$ for any $t'\geq t$. The collection of partitions $(\Pi^{0}(t))_{t\geq 0}$ will be a coalescent process. To define these processes and go from an exchangeable distinguished partition to another coarser partition, we have to introduce the coagulation operator.
\subsection{Coagulation operator and distinguished coalescents} 
To define the distinguished exchangeable coalescents, we need to define an operator on the space of distinguished partitions.
\begin{Def} 
Let $\pi, \pi'  \in \mathcal{P}^{0}_{\infty}$, the partition $Coag(\pi,\pi')$ is defined by $Coag(\pi,\pi')_{i}=\pi''_{i}$ where $\pi''_{i}=\bigcup_{j\in \pi'_{i}}\pi_{j}$. The partition $Coag(\pi,\pi')$ is exactly the one obtained by coagulating blocks of $\pi$ according to blocks of $\pi'$. 
\end{Def}
We denote by $0_{\cro{\infty}}$ the partition into singletons $\{\{0\},\{1\},...\}$, and by $1_{\cro{\infty}}$ the trivial partition $\{\mathbb{Z}^{+},\emptyset,...\}$. Plainly, for all $n\geq 0$, $Coag(\pi,\pi')_{|\cro{n}}=Coag(\pi_{|\cro{n}},\pi'_{|\cro{n}})$ and for all $\pi \in \mathcal{P}_{\infty}^{0}, Coag(\pi,0_{\cro{\infty}})=\pi$ and $Coag(\pi,1_{\cro{\infty}})=1_{\cro{\infty}}$. Note that however, we do not have $Coag(\pi,\pi^{'})_{|K}=Coag(\pi_{|K},\pi^{'}_{|K})$ for $K\subset \mathbb{N}$ in general.
\begin{prop} Let $\pi$, $\pi'$ be two independent exchangeable distinguished partitions. The distinguished partition $Coag(\pi,\pi')$ is also exchangeable. 
\end{prop}
\textit{Proof.} See the proof of Lemma 4.3 in \cite{coursbertoin}. 
\\
The coagulation operator allows us to define distinguished coalescents which are Markovian processes valued in distinguished partitions of $\mathbb{N}$.
\begin{Def} A Markov process $\Pi^{0}$ with values in $\mathcal{P}^{0}_{\infty}$ is called a distinguished coalescent if its semi-group is given as follows: the conditional law of $\Pi^{0}(t+t')$ given $\Pi^{0}(t)=\pi^{0}$ is the law of $Coag(\pi^{0},\pi')$ where $\pi'$ is some exchangeable distinguished partition (whose law only depends on $t'$). A distinguished coalescent is called standard if $\Pi^{0}(0)=0_{\textbf{[}\infty\textbf{]}}$.
\end{Def}
The properties of the coagulation operator, see \cite{coursbertoin}, imply that a distinguished coalescent $\Pi^{0}$ fulfills the Feller property. Therefore the process has a càdlàg version and is strong Markovian. Plainly, the random partition $\Pi^{0}_{|\mathbb{N}}(t)$ is exchangeable for all $t\geq 0$. However, we stress that in general the process $(\Pi^{0}_{|\mathbb{N}}(t),t\geq 0)$ is not an exchangeable coalescent and not even Markovian. We will give an example in Section 4.2. \\
For every $n\geq 1$, the restriction $\Pi^{0}_{|\cro{n}}$ is a continuous time Markov chain with a semi-group given by the operator $Coag$. Let $\pi\in \mathcal{P}^{0}_{n}\setminus \{0_{\cro{n}}\}$, we denote by $q_{\pi}$, the jump rate of $\Pi^{0}_{|\cro{n}}$ from $0_{\cro{n}}$ to $\pi$: $$q_{\pi}:=\lim_{t\rightarrow 0+} \frac{1}{t}\mathbb{P}_{0_{\cro{n}}}[\Pi^{0}_{|\cro{n}}(t)=\pi].$$
An easy adaptation of the proof of Proposition 4.4 in \cite{coursbertoin} gives the existence and uniqueness of the distinguished coagulation measure:
\begin{Def} \label{coagrate} 
The distinguished coagulation measure of $\Pi^{0}$ is the unique measure $\mu^{0}$ on $\mathcal{P}^{0}_{\infty}$, such that $\mu^{0}(\{0_{\cro{\infty}}\})=0$ and
$$\mu^{0}(\mathcal{P}^{0}_{\infty,\pi})=q_{\pi}$$
for every $n\in \mathbb{Z}_{+}$ and every partition $\pi \in \mathcal{P}^{0}_{n}$. 
\\
Moreover, the measure $\mu^{0}$ fulfills
\begin{center}
$\mu^{0}(\pi \in \mathcal{P}^{0}_{\infty}: \pi_{|\cro{n}}\neq 0_{\cro{n}})<\infty$ and $\mu^{0}$ is exchangeable.
\end{center}
Conversely, any measure fulfilling the previous conditions will be called a distinguished coagulation measure.
\end{Def}
Let $\mu^{0}$ be a distinguished coagulation measure, we construct explicitly a distinguished coalescent process with coagulation measure (in the sense of Definition \ref{coagrate}) $\mu^{0}$. Let $N$ be a Poisson measure with intensity $dt\otimes \mu^{0}(d\pi)$. Let $N_{b}$ be the image of $N$ by the map $(t,\pi)\mapsto (t, \pi_{|\cro{b}})$. Its intensity, denoted by $\mu^{0}_{b}$, is the image of $\mu^{0}$ by the previous map. We denote by $(t_{i},\pi_{i})$ the atoms of $N_{b}$ and define a process $(\Pi^{b}(t), t\geq 0)$ by the following recursion:
\begin{center} For all $0\leq t<t_{1}$  $\Pi^{b}(t)=0_{\cro{b}}$ and if $t_{i}\leq t<t_{i+1}$, $\Pi^{b}(t)=Coag(\Pi^{b}(t_{i-1}),\pi^{i})$, with $t_{0}=0$.
\end{center}
\begin{prop} \label{construction} The sequence of random partitions $(\Pi^{b}(t), b\in \mathbb{N})$ is compatible, which means that for all $a\leq b$, $\Pi^{b}_{|\cro{a}}=\Pi^{a}$. The unique process $(\Pi^{0}(t), t\geq 0)$ such that $\Pi^{0}_{|\cro{b}}(t)=\Pi^{b}(t)$, defined by $\Pi^{0}_{i}(t)=\bigcup_{b\geq 1}\Pi^{b}_{i}(t)$ is a distinguished coalescent with coagulation measure $\mu^{0}$. 
\end{prop}
\textit{Proof.} Same arguments as those of Proposition 4.5 in \cite{coursbertoin} apply.
$\square$
\begin{Example} We denote by $K(i,j)$ the simple distinguished partition where $i$ and $j$ are in the same block and all the other blocks are singletons. Let $c_{0}$, $c_{1}$ be two non-negative real numbers. The measure $\mu^{0}=c_{0}\mu_{0}^{K}+c_{1}\mu_{1}^{K}$, where $\mu_{0}^{K}:=\sum_{1\leq i} \delta_{K(0,i)}$ and $\mu_{1}^{K}:=\sum_{1\leq i<j} \delta_{K(i,j)}$, is a distinguished coagulation measure. The process obtained is called the Kingman's distinguished coalescent with rates $(c_{0},c_{1})$.
\end{Example}
Indeed, the measure $\mu^{0}$ defined as above is plainly a distinguished coagulation measure. The Poissonian construction explains the dynamics of this process. At a constant rate $c_{0}$, a block not containing $0$ merges with $\Pi^{0}_{0}$ that is a singular coagulation with the distinguished block. At a constant rate $c_{1}$, two blocks not containing $0$ merge into one, that is the classic binary coagulation of Kingman's coalescent.  
\\
\subsection{Characterization in law of the distinguished coalescents}
The next theorem is one of the main results of this work, it claims that the law of a distinguished coalescent is characterized by two non-negative real numbers $c_{0}$, $c_{1}$ and a measure $\nu^{0}$ on $\mathcal{P}^{0}_{\textbf{m}}$. It should be viewed as an extension of Theorem \ref{repres} to certain infinite measures. Recall that $\rho^{0}_{s}$ denotes the law of an exchangeable distinguished $s$-paint-box for $s\in \mathcal{P}^{0}_{ \textbf{m}}$.
\begin{thm} \label{maintheo} Recalling Definition \ref{coagrate}, let $\mu^{0}$ be a distinguished coagulation measure. There exist two unique real numbers $c_{0}, c_{1}$ and a unique measure $\nu^{0}$ on $\mathcal{P}^{0}_{\textbf{m}}$ which fulfills:
\begin{center}
$\nu^{0}(0)=0$ and $\int_{\mathcal{P}^{0}_{\textbf{m}}} (s_{0}+\sum_{i=1}^{\infty}s^{2}_{i})\nu^{0}(ds)< \infty$ 
\end{center}
such that $$\mu^{0}=c_{0}\mu_{0}^{K}+c_{1}\mu_{1}^{K}+\rho^{0}_{\nu^{0}}$$
where $$\rho^{0}_{\nu^{0}}(.):=\int_{s\in \mathcal{P}^{0}_{\textbf{m}}}\rho^{0}_{s}(.)\nu^{0}(ds).$$
Conversely, let $c_{0}, c_{1}$ and $\nu^{0}$ be two real numbers and a measure on $\mathcal{P}^{0}_{\textbf{m}}$ verifying the previous conditions, there exists a unique (in law) distinguished coalescent with $\mu^{0}=c_{0}\mu_{0}^{K}+c_{1}\mu_{1}^{K}+\rho^{0}_{\nu^{0}}$.
\end{thm}
When $\nu^{0}$ is carried on $\{s\in \mathcal{P}^{0}_{\textbf{m}}; s_{0}=0\}$ (which can be identified as $\mathcal{P}_{\textbf{m}}$), the block containing $0$ is reduced to the singleton $\{0\}$ (we distinguish no block) and considering the restriction to $\mathbb{N}$, we recover the characterization of exchangeable coalescents (also called $\Xi$-coalescents) by Schweinsberg in \cite{schweinsberg}. 
\\
\\
\textit{Proof.} Arguments used to prove this theorem are adapted from those of Theorem 4.2 in Chapter 4 of \cite{coursbertoin}. Nevertheless, we give details to highlight the fact that the condition on $\nu^{0}$ differs from that of Theorem 4.2 in \cite{coursbertoin}. We denote by $\mu^{0}_{n}$ the restriction of $\mu^{0}$ to $\{\pi \in \mathcal{P}_{\infty}^{0}; \pi_{|\cro{n}} \neq 0_{\cro{n}}\}$. The measure $\mu^{0}_{n}=1_{\{\pi \in \mathcal{P}_{\infty}^{0}; \pi_{|\textbf{[}n\textbf{]}}\neq 0_{\cro{n}}\}}\mu^{0}$ has a finite mass and is invariant under the action of permutations $\sigma$ that coincide with the identity on $\cro{n}$.
We define the $n$-shift on distinguished partitions by the map $\pi \to \pi'$ defined by : 
\begin{align*}
\forall i,j \geq 1: \quad &i\underset{\pi'}{\sim} j \Longleftrightarrow i+n \underset{\pi}{\sim} j+n \\
\forall j \geq 1: \quad &0\underset{\pi'}{\sim} j \Longleftrightarrow 0 \underset{\pi}{\sim} j+n.
\end{align*}
The image of $\mu^{0}_{n}$ by the $n$-shift, denoted by $\bar{\mu}^{0}_{n}$, is invariant under the action of permutations $\sigma$ of $\mathbb{Z}_{+}$ such that $\sigma(0)=0$. By the Kingman's correspondence (Theorem \ref{repres})
$$\bar{\mu}^{0}_{n}(d\pi)=\int_{\mathcal{P}_{\textbf{m}}^{0}}\rho_{s}^{0}(d\pi)\bar{\mu}^{0}_{n}(|\pi|^{\downarrow}\in ds).$$ 
Moreover, $\bar{\mu}^{0}_{n}$ almost every partition has asymptotic frequencies. The shift does not affect asymptotic frequencies and then $\mu^{0}_{n}$ a.e partition has asymptotic frequencies. The measure $\mu^{0}$ is the increasing limit of the $\mu^{0}_{n}$, we deduce that $\mu^{0}$ almost every partition has asymptotic frequencies. \\
By the distinguished paint-box representation of $\bar{\mu}^{0}_{n}$, we get for all $s\in \mathcal{P}_{\textbf{m}}^{0}\setminus \{0\}$
\begin{center}
$\mu^{0}_{n}(n+1\sim n+2 \not \sim 0$ or $0\sim n+1)| |\pi|^{\downarrow}=s)=s_{0}+\sum_{k=1}^{\infty}s^{2}_{k}$.
\end{center}
Let $\nu^{0}_{n}$ be the image of $\mu^{0}_{n}$ by the map $\pi \mapsto |\pi|^{\downarrow}$
$$\nu^{0}_{n}(ds)=\mu_{n}^{0}(|\pi|^{\downarrow}\in ds).$$
We stress that $\nu^{0}_{n}(ds)=\bar{\mu}^{0}_{n}(|\pi|^{\downarrow}\in ds)$ because the $n$-shift has no impact on asymptotic frequencies.
\\
We have : $\mu^{0}_{n}(n+1 \sim n+2 \not \sim 0$ or $0\sim n+1)\geq \int_{\mathcal{P}^{0}_{\textbf{m}}}(s_{0}+\sum_{i=1}^{\infty}s^{2}_{i})\nu_{n}^{0}(ds)$.
\\
Moreover, $\mu^{0}_{n}(n+1 \sim n+2 \not \sim 0$ or $0\sim n+1)\leq \mu^{0}(n+1 \sim n+2 \not \sim 0$ or $0\sim n+1)$.
\\
By exchangeability of $\mu^{0}$
\begin{center} 
$\mu^{0}(n+1 \sim n+2 \not \sim 0$ or $0\sim n+1)=\mu^{0}(1\sim 2 \not \sim 0$ or $0 \sim 1)\leq\mu^{0}(\pi_{|\cro{2}}\neq 0_{\cro{2}})<\infty$.
\end{center}
We deduce that the finite measures $\nu^{0}_{n}$ increase as $n\uparrow \infty$ to the measure $\nu^{0}:=\mu^{0}(|\pi|^{\downarrow}\in ds)$ and so $$\underset{n\rightarrow \infty}{\lim}\int_{\mathcal{P}^{0}_{\textbf{m}}}(s_{0}+\sum_{i=1}^{\infty}s^{2}_{i})\nu_{n}^{0}(ds)=\int_{\mathcal{P}^{0}_{\textbf{m}}}(s_{0}+\sum_{i=1}^{\infty}s^{2}_{i})\nu^{0}(ds)\leq \mu^{0}(\pi_{|\cro{2}}\neq 0_{\cro{2}})<\infty .$$
Let $k \in \mathbb{N}$ and $\pi^{[k]} \in \mathcal{P}_{k}\setminus \{0_{\textbf{[}k\textbf{]}}\}$. The sequence of events $\pi_{|\{k+1,...,k+n\}}\neq 0_{\{k+1,...,k+n\}}$ is increasing, then we have
$$\mu^{0}(\pi_{|\textbf{[}k\textbf{]}}=\pi^{[k]}, |\pi|^{\downarrow}\neq 0)= \underset{n\rightarrow \infty}{\lim} \mu^{0}(\pi_{|\textbf{[}k\textbf{]}}=\pi^{[k]},|\pi|^{\downarrow}\neq 0, \pi_{|\{k+1,...,k+n\}}\neq 0_{\{k+1,...,k+n\}}).$$
By an obvious permutation we get
$$\mu^{0}(\pi_{|\textbf{[}k\textbf{]}}=\pi^{[k]},|\pi|^{\downarrow}\neq 0, \pi_{|\{k+1,...,k+n\}}\neq 0_{\{k+1,...,k+n\}})=\bar{\mu}_{n}^{0}(\pi_{|\textbf{[}k\textbf{]}}=\pi^{[k]},|\pi|^{\downarrow}\neq 0).$$
Thus, using the distinguished paint-box representation of $\bar{\mu}_{n}^{0}$, we deduce that 
$$\bar{\mu}_{n}^{0}(\pi_{|\textbf{[}k\textbf{]}}=\pi^{[k]},|\pi|^{\downarrow}\neq 0)=\int_{\mathcal{P}_{\textbf{m}}^{0}}\rho_{s}^{0}(\pi_{|\textbf{[}k\textbf{]}}=\pi^{[k]},|\pi|^{\downarrow}\neq 0)\nu_{n}^{0}(ds) \underset{n\rightarrow \infty}{\longrightarrow} \int_{\mathcal{P}_{\textbf{m}}^{0}}\rho_{s}^{0}(\pi_{|\textbf{[}k\textbf{]}}=\pi^{[k]},|\pi|^{\downarrow}\neq 0)\nu^{0}(ds).$$ 
As $k$ is arbitrary, we get 
$$1_{\{|\pi|^{\downarrow}\neq 0 \}} \mu^{0}(d\pi)=\int_{\mathcal{P}_{\textbf{m}}^{0}}\rho_{s}^{0}(d\pi)\nu^{0}(ds).$$
It remains to study $1_{\{|\pi|^{\downarrow}=0 \}} \mu^{0}(d\pi)$. Consider now $\tilde{\mu}^{0}(d\pi):=1_{\{0\sim 1, |\pi|^{\downarrow}=0\}}\mu^{0}(d\pi)$ which has finite mass (because $\mu^{0}(0\sim 1)<\infty$). We want to show that $\tilde{\mu}^{0}(d\pi)$ is proportional to $\delta_{K(0,1)}$, where $K(0,1)$ is the simple partition with $0\sim 1$. Let $\tilde{\mu}^{0}_{2}(d\pi)$ be the image of $\tilde{\mu}^{0}$ by the $2$-shift. The measure $\tilde{\mu}^{0}_{2}(d\pi)$ is supported by $\{ \pi \in \mathcal{P}_{\infty}; |\pi|^{\downarrow}=0 \}$ and is exchangeable with finite mass. By the distinguished paint-box construction, the only exchangeable partition with asymptotic frequencies $|\pi|^{\downarrow}=0$ is the partition into singletons $0_{\textbf{[}\infty\textbf{]}}$. Therefore, $\tilde{\mu}^{0}_{2}(d\pi)=c_{0}\delta_{ \{\{0\},\{1\},....\} }$. We deduce that for $\tilde{\mu}^{0}_{2}$-almost every $\pi'$, $\forall i\neq j$, $i\not \underset{\pi'}{\sim} j$.  From the definition of the $2$-shift, we get that for $\tilde{\mu}^{0}$-almost every $\pi$ 
\begin{center} $\forall i,j\geq 1$; $i \neq j$, $i+2 \not \underset{\pi}{\sim} j+2$, and $\forall j\geq 1$, $0 \not \underset{\pi}{\sim} j+2$. \end{center}
It implies that we have to consider only three possibilities \begin{center}
$\tilde{\mu}^{0}(d\pi)=c_{0}\delta_{K(0,1)}$ or $\tilde{\mu}^{0}(0\sim 1\sim 2)>0$ or $\tilde{\mu}^{0}(2\sim k)>0$ for some $k\geq 3$.\end{center}
If $\tilde{\mu}^{0}(2\sim k)>0$ for some $k\geq 3$, we get by exchangeability that $\tilde{\mu}^{0}(2\sim k)=\tilde{\mu}^{0}(2\sim 3)>0$. Moreover, the collection of sets $\{\{2\sim n\}, n\geq 3\}$ is such that the intersection of two or more sets has a zero measure $\tilde{\mu}^{0}$ and then $\tilde{\mu}^{0}(\cup_{n\geq 3} \{2\sim n\})=\sum_{n\geq 3}\tilde{\mu}^{0}(2\sim n)\leq \mu^{0}(0\sim 1)$. It follows that $\mu^{0}(0\sim 1)=\infty$. This is a contradiction because $\mu^{0}(0\sim 1)<\infty$.
If $\tilde{\mu}^{0}(0\sim 1\sim 2)>0$, then by exchangeability for all $n\geq 2,$ $\tilde{\mu}^{0}(0\sim 1\sim n)=c_{0}>0$ and by the same arguments, the same contradiction appears. We deduce that $\tilde{\mu}^{0}(d\pi)$ is equal to $c_{0}\delta_{K(0,1)}$. By exchangeability, we have $1_{\{0\sim i, |\pi|^{\downarrow}=0\}}\mu^{0}(d\pi)=c_{0}\delta_{K(0,i)}$. The measure $\mu^{0}1_{\{|\pi|^{\downarrow}=0, \ \pi_{0}\neq \{0\}\}}$ is carried on simple partition $\pi$ such that $\pi_{0}$ is not singleton, moreover the collection of sets $\{\{0\sim i\}, i\geq 1\}$ is such that the intersection of two sets has a zero measure. Therefore, we have $$1_{\{|\pi|^{\downarrow}=0,\ \pi_{0}\neq \{0\} \}}\mu^{0}(d\pi)=c_{0}\sum_{i\geq 1}\delta_{K(0,i)}.$$ 
The restriction of $\mu^{0}$ to $\{\pi \in \mathcal{P}^{0}_{\infty}; \pi_{0}= \{0\}\}$ can be viewed as an exchangeable measure on $\mathcal{P}_{\infty}$, the argument to conclude is then the same as in \cite{coursbertoin} on p184.
$\square$
\begin{rem}
Denoting by $(D_{0}(t))_{t\geq 0}:=(1-\sum_{i=0}^{\infty}|\Pi^{0}_{i}(t)|)_{t\geq 0}$ the process of dust of $\Pi^{0}$. The arguments of \cite{Pitman2} or \cite{coursbertoin} allow to show that for all $t>0$, the random partition $\Pi^{0}(t)$ has improper asymptotic frequencies with a strictly positive probability if and only if $c_{1}=0$ and $\int_{\mathcal{P}^{0}_{\textbf{m}}} (1-\delta)\nu^{0}(ds)<\infty$ where $\delta=1-\sum_{i=0}^{\infty}s_{i}$. \\
In that case, the process $(\xi^{0}(t))_{t\geq 0}:=(-ln(D_{0}(t)))_{t\geq 0}$ is a subordinator with Laplace exponent $$\phi^{0}(q)=c_{0}q+\int_{\mathcal{P}^{0}_{\textbf{m}}}(1-\delta^{q})\nu^{0}(ds).$$ 
Note that the drift coefficient $c_{0}$ may be positive, which contrasts with the result of Pitman \cite{Pitman2}.
\end{rem} 

\section{The simple distinguished exchangeable coalescents : $M$-coalescents}
In this section, we focus on simple distinguished coalescent for which the coagulation measure $\mu^{0}$ is carried by the set of simple distinguished partitions. We call them, hereafter, $M$-coalescents. These processes are the analogue of $\Lambda$-coalescents for exchangeable distinguished coalescents. Historically, the $\Lambda$-coalescent is the first exchangeable coalescent with multiple collisions to have been defined, see Pitman, \cite{Pitman2} and Sagitov, \cite{Sagitov}. We begin by recalling some basic facts about $\Lambda$-coalescents.
\subsection{$\Lambda$-coalescents}
A $\Lambda$-coalescent (also called simple exchangeable coalescent) is a process taking values in the partitions of $\mathbb{N}$ describing the genealogy of an infinite haploid population, labelled by $\mathbb{N}$ where two or more ancestral lineages merging cannot occur simultaneously. We stress that in these coalescent processes, each individual has an ancestor in the population. Immigration phenomenon is not taken into account and no block is distinguished. A simple exchangeable coalescent is a Markovian process $(\Pi(t),t\geq 0)$ on the space of partitions of $\mathbb{N}$ satisfying :\\
\begin{itemize}
\item[i)] If $n\in \mathbb{N}$, then the restriction $(\Pi_{|[n]}(t))$ is a continuous-time Markov chain valued in $\mathcal{P}_{n}$;\\
\item[ii)] For each $n$, $(\Pi_{|[n]}(t))$ evolves by exchangeable merging of blocks : $\Pi_{|[n]}(t)=Coag(\Pi_{|[n]}(t-),\pi')$ where $\pi'$ is an independent simple exchangeable partition;\\ 
\end{itemize}
By Theorem 1 in Pitman \cite{Pitman2}, or Sagitov \cite{Sagitov}, we know that any simple exchangeable coalescent is characterized in law by a finite measure $\Lambda$ on $[0,1]$. The dynamics of $\Pi$ can be described as follows: whenever $\Pi_{|[n]}(t)$ is a partition with $b$ blocks, the rate at which a $k$-tuple of its blocks merges is
$$\lambda_{b,k}=\int_{0}^{1}x^{k-2}(1-x)^{b-k}\Lambda(dx).$$
When $\Lambda$ is the Dirac at $0$, we recover the Kingman's coalescent. When $\Lambda(\{0\})=0$, the $\Lambda$-coalescent can be constructed via a Poisson Point Process on $\mathbb{R}_{+}\times [0,1]$ with intensity $dt\otimes \nu(dx)$, where $\nu(dx)=x^{-2}\Lambda(dx)$ :$$N=\sum_{i\in \mathbb{N}}\delta_{(t_{i},x_{i})}.$$
The atoms of $N$ encode the evolution of the coalescent $\Pi$: At time $t-$, flip a coin with probability of "heads" $x$ for each block. All blocks flipping "heads" are merged immediately. We can, also, construct a simple exchangeable partition, $\pi'$ where the non trivial block is constituted by indices of "heads". Thus, we get $\Pi(t)$ by $Coag(\Pi(t-),\pi')$. In order to make this construction rigorous, one first considers the restrictions $(\Pi_{|[n]}(t))$ as in Proposition \ref{construction}, since the measure $\nu(dx):=x^{-2}\Lambda(dx)$ can have an infinite mass.  
\subsection{$M$-coalescents}
The distinguished exchangeable coalescents such that when a coagulation occurs all the blocks involved merge as a single block are called $M$-coalescents. We specify their laws by two finite measures on $[0,1]$, and study their generators in the same fashion as those of $\Lambda$-coalescents.
\begin{Def} When a distinguished coagulation measure $\mu^{0}$ is carried by the set of simple distinguished partitions (with only one block non empty nor singleton), the distinguished coalescent $\Pi^{0}$ is said to be simple. Define the following restricted measures:
\begin{center} $\nu_{0}=\nu^{0}1_{\{s\in \mathcal{P}^{0}_{\textbf{m}};s=(s_{0},0,...)\}}$ and $\nu_{1}=\nu^{0}1_{\{s\in \mathcal{P}^{0}_{\textbf{m}};s=(0,s_{1},0,...)\}}$. 
\end{center}
We can write $\nu^{0}=\nu_{0}+\nu_{1}$. By a slight abuse of notation, $\nu_{0}$ and $\nu_{1}$ can be viewed as two measures on $[0,1]$ such that $\int_{0}^{1}s_{0}\nu_{0}(ds_{0})<\infty$ and $\int_{0}^{1}s_{1}^{2}\nu_{1}(ds_{1})<\infty$, and Theorem \ref{maintheo} yields
$$\mu^{0}=c_{0}\mu_{0}^{K}+\rho^{0}_{\nu_{0}}+c_{1}\mu_{1}^{K}+\rho_{\nu_{1}}.$$
We define the finite measures $\Lambda_{0}(dx):=x\nu_{0}(dx)+c_{0}\delta_{0}$, and $\Lambda_{1}(dx):=x^{2}\nu_{1}(dx)+c_{1}\delta_{0}$. The law of a simple distinguished coalescent is then characterized by $M=(\Lambda_{0},\Lambda_{1})$, and we call this subclass the $M$-coalescents.
\end{Def}
As already mentioned in Section 3, in the most cases, the restriction to $\mathbb{N}$ of a distinguished coalescent is not Markovian. Let $\Pi^{0}$ be a $M$-coalescent with for instance $\Lambda_{0}=\delta_{0}$ and $\Lambda_{1}(dx)=dx$ (the Lebesgue measure). To locate the distinguished block in $\Pi^{0}_{|\mathbb{N}}$, we may locate a binary coagulation before time $t$ (all other mergers involve an infinite number of blocks). The restricted process $\Pi^{0}_{|\mathbb{N}}$ is then not Markovian.
\\
\\
The explicit Poissonian construction of Proposition \ref{construction} can now be interpreted in the same way as the one of $\Lambda$-coalescent, recalled in Section 4.1. When $\Lambda_{0}(\{0\})=\Lambda_{1}(\{0\})=0$, the $M$-coalescent associated can be constructed via two Poisson Point Processes $N_{0}$ and $N_{1}$ on $\mathbb{R}_{+}\times (0,1]$ with intensities $dt\otimes \nu_{0}(dx)$ and $dt\otimes \nu_{1}(dx)$, where $\nu_{0}(dx)=x^{-1}\Lambda_{0}(dx)$ and $\nu_{1}(dx)=x^{-2}\Lambda_{1}(dx)$: 
\begin{itemize}
\item At an atom $(t_{i},x_{i})$ of $N_{1}$, flip a coin with probability of "heads" $x_{i}$ for each block not containing $0$. All blocks flipping "heads" are merged immediately in one block as in the Proposition \ref{construction}. 
\item At an atom $(t_{i},x_{i})$ of $N_{0}$, flip a coin with probability of "heads" $x_{i}$ for each block not containing $0$. All blocks flipping "heads" coagulate immediately with the distinguished block. 
\end{itemize}
This construction is exactly the one obtained when we coagulate the partition at $t-$ with a simple exchangeable distinguished partition $\pi'$ where the non trivial block is constituted by indexes of "heads". Thus, we construct the $M$-coalescent in the same way that $\Lambda$-coalescent in Section 4.1.
\\
\\
We investigate jump rates of a $M$-coalescent $(\Pi^{0}(t))_{t\geq 0}$. Thanks to the simple distinguished paint-box structure, we compute explicitly the jump rates of the restriction of $\Pi^{0}$. Let $\pi\in \mathcal{P}^{0}_{n}$ be simple, $q_{\pi}=\mu^{0}(\mathcal{P}^{0}_{\infty},\pi)$:
\begin{itemize}
\item For every $2\leq k\leq n$, if $\pi$ has one block not containing $0$ with $k$ elements, then $$q_{\pi}=\lambda_{n,k}:=\int_{0}^{1}x^{k-2}(1-x)^{n-k}\Lambda_{1}(dx).$$
\item For every $1\leq k\leq n$, if the distinguished block of $\pi$ has $k+1$ elements (counting $0$), then $$q_{\pi}=r_{n,k}:=\int_{0}^{1}y^{k-1}(1-y)^{n-k}\Lambda_{0}(dy).$$
\end{itemize}
Let $\pi \in \mathcal{P}^{0}_{p}$ with $b$ blocks without $0$ and $F$ any function defined on $\mathcal{P}^{0}_{p}$. The generator of $\Pi^{0}_{|\cro{p}}$ is
$$\mathcal{L}^{*}F(\pi)=\sum_{I\subset \{1,...,b\}, |I|\geq 2}\lambda_{b,|I|}(F(c_{I}\pi)-F(\pi))+\sum_{J\subset \{1,...,b\}, |J|\geq 1}r_{b,|J|}(F(c^{J}\pi)-F(\pi)),$$
with $c_{I}\pi=Coag(\pi,\{\{0\},\{1\},...,\{I\},...\})$ and $c^{J}\pi=Coag(\pi,\{\{0\}\cup\{J\},\{.\},...,\{.\}\}).$
\subsection{Coming down from infinity for $M$-coalescents}
Let $\Lambda$ be a finite measure on $[0,1]$ and $\Pi$ be a $\Lambda$-coalescent. Pitman \cite{Pitman2} showed that if $\Lambda(\{1\})=0$, only the following two types of behavior are possible, either $\mathbb{P}[$For all $t>0$, $\Pi(t)$ has infinitely many blocks$]=1$, or $\mathbb{P}[$For all $t>0$, $\Pi(t)$ has only finitely many blocks$]=1$. In the second case, the process $\Pi$ is said to come down from infinity. For instance, Kingman's coalescent comes down from infinity, while if $\Lambda(dx)=dx$, then the corresponding $\Lambda$-coalescent (called Bolthausen-Sznitman coalescent) does not come down from infinity. A necessary and sufficient condition for a $\Lambda$-coalescent to come down from infinity was given by Schweinsberg in \cite{CDI}. Define $$\phi(n)=\sum_{k=2}^{n}(k-1)C^{k}_{n}\lambda_{n,k}$$ with $\lambda_{n,k}=\int_{0}^{1}x^{k-2}(1-x)^{n-k}\Lambda(dx)$. The $\Lambda$-coalescent comes down from infinity if and only if $\sum_{n=2}^{\infty}\frac{1}{\phi(n)}<\infty$.\\
Define $\psi_{\Lambda}(q):=\int_{[0,1]}(e^{-qx}-1+qx)x^{-2}\Lambda(dx)$, Bertoin and Le Gall observed, in \cite{LGB3} the following equivalence
$$\sum_{n=2}^{\infty}\frac{1}{\phi(n)}<\infty \Longleftrightarrow \int_{a}^{\infty}\frac{dq}{\psi_{\Lambda}(q)}<\infty$$
where the right-hand side holds for some $a>0$ (and then necessary for all). This equivalence is explained in a probabilistic way by Berestycki et al. in \cite{Beres}. 
\\
\\
As for the $\Lambda$-coalescent, if the $M$-coalescent comes down from infinity, it does immediately:
\begin{prop} \label{coming}
Let $(\Pi^{0}(t))_{t\geq 0}$ be a $M$-coalescent, with $\Lambda_{0}$ and $\Lambda_{1}$ without mass at $1$. We denote by $T$ its time of coming down from infinity: $T=\inf \{t>0, \#\Pi(t)<\infty\}$. We have a.s $T=0$ or $T=\infty$. 
\end{prop}
\textit{Proof.} See the proof of Lemma 31 in \cite{schweinsberg}. $\square$
\\
We stress that when $\Lambda_{0}+\Lambda_{1}$ has a mass at $1$, the $M$-coalescent comes down from infinity. Indeed, by the Poisson construction, in an exponential time $\tau$ of parameter $(\Lambda_{0}+\Lambda_{1})(\{1\})$, the Poisson measure $N$ has an atom $\pi$ such that $\pi_{0}=\mathbb{Z}_{+}$ or $\pi_{1}=\mathbb{N}$. Thus for large $t$, the process $\Pi^{0}(t)$ has just one block.
\\
\\
It remains to focus on the case where $\Lambda_{0}+\Lambda_{1}$ has no mass at $1$. Intuitively, when the genuine $\Lambda_{1}$-coalescent comes down, all blocks merged in one in an almost surely finite time. On the one hand, we can think that the $(\Lambda_{0},\Lambda_{1})$-coalescent has more jumps and coagulates all its blocks faster. On the other hand, the perturbation due to the coagulation with the distinguished block on the general term of the sum, studied initially by Schweinsberg in \cite{CDI}, is not sufficient to induce its convergence and so the coming down. The $(\Lambda_{0}, \Lambda_{1})$-coalescent comes down from infinity if and only if the $\Lambda_{1}$-coalescent comes down: 
\begin{thm} \label{cdi} The $(\Lambda_{0},\Lambda_{1})$-coalescent comes down from infinity if and only if $\sum_{n=2}^{\infty}\frac{1}{\phi_{1}(n)}<\infty$ where 
$\phi_{1}(n)=\sum_{k=2}^{n}(k-1)C^{k}_{n}\lambda_{n,k}$ and $\lambda_{n,k}$ as in Section 4.2.
\end{thm}
The proof requires rather technical arguments and is given in the Section 6. 
\section{M-coalescents and generalized Fleming-Viot processes with immigration}
In the final section, we are interested in a population model which has exactly a genealogy given by a $M$-coalescent. A powerful method to study simultaneously the population model and its genealogy is to define some stochastic flows as Bertoin and Le Gall in \cite{LGB1}. A process valued in the space of \textit{probability measures} on $[0,1]$, $(Z^{0}_{t},t\geq 0)$ is embedded in the flow. The atoms of the random probability $Z^{0}_{t}$ represent the current types frequencies in the population at time $t$. Moreover, $Z^{0}_{t}$ has a distinguished atom at $0$ representing the fraction of immigrants in the population. This process will be called $M$-generalized Fleming-Viot processes with immigration. Following \cite{LGB1}, we begin by establishing a correspondence between some stochastic flows and $M$-coalescents. 
\subsection{Stochastic flows of distinguished bridges}
By assumption, at any time the families describing the population form a distinguished exchangeable partition. Theorem \ref{repres} ensures that it has a distinguished paint-box structure. We have to study some random functions called distinguished bridges.
\subsubsection{Distinguished bridges and exchangeable distinguished partitions}
Considering the underlying law on $[0,1]$ associated with a $s$-distinguished paint-box (see definition \ref{paintbox}), we introduce the distinguished bridges defined by
$$b_{s}(r)=s_{0}+\sum_{i=1}^{\infty}s_{i}1_{\{V_{i}\leq r\}}+\delta r$$ 
where $s$ is a distinguished mass-partition and $(V_{i})_{i\geq 1}$ a sequence of independent uniform variables. 
Let $U_{0}=0$ and $(U_{i})_{i\geq 1}$ be an independent sequence of i.i.d uniform variables. The partition given by $i\sim j$ iff $b_{s}^{-1}(U_{i})=b_{s}^{-1}(U_{j})$ is exactly the $s$-distinguished paint-box. When $s_{0}=0$, the bridge encodes a paint-box partition with no distinguished block.
\begin{figure}[!h]
\centering 
\includegraphics[height= .25 \textheight ] {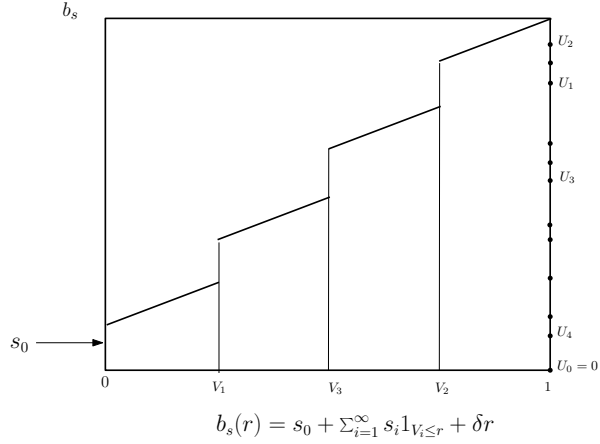}
\caption{Distinguished bridge}
\end{figure}
\\
Focussing on $M$-coalescents, we only need to focus on two types of distinguished bridges: 
\begin{itemize}
\item Bridges with distinguished mass-partition $(0,x,0,...)$: $b_{0,x}(r)=x1_{V\leq r}+r(1-x)$. 
\item Bridges with distinguished mass-partition $(x,0,0,...)$ : $b_{x,0}(r)=x+r(1-x)$.
\end{itemize}
To be concise, we shall work directly with distinguished bridges of the form $$b_{y,x}(r)=y+x1_{V\leq r}+r(1-x-y).$$ 
The following lemma relates the composition of distinguished bridges to the coagulation of simple distinguished partitions.
\begin{lem} \label{key} Let $x,x',y,y' \in [0,1]$  with $x+y,x'+y' \leq 1$ and $b_{y,x}$, $b_{y',x'}$ two independent distinguished bridges. \\
Let $\pi$ : $i \sim j$ iff $b_{y,x}^{-1}(U_{i})=b_{y,x}^{-1}(U_{j})$. We stress that $\pi$ has at most two non-trivial blocks. 
\begin{itemize}
\item [(i)]  For $i\geq 1$, we set $U'_{i}:=b_{y,x}^{-1}(U_{j})$, $\forall j\in \pi_{i}$. The variables $(U'_{i})_{i\geq 1}$ are i.i.d uniform, independent of $\pi$.
\item [(ii)] Let $\pi'$ be the partition constructed from $b_{y',x'}$ and $(U'_{i})_{i\geq 1}$, we denote by $\pi^{0}$ the partition such that \\
$i \sim j$ iff $b_{y,x}^{-1}\circ b_{y',x'}^{-1}(U_{i})=b_{y,x}^{-1}\circ b_{y',x'}^{-1}(U_{j})$. We have the identity $\pi^{0}=Coag(\pi,\pi')$. 
\end{itemize}
\end{lem}
\textit{Proof}. The proof is an easy adaptation of Lemma 4.8 of \cite{coursbertoin}.$\square$
\\
\\
This lemma is the key observation in order to associate a stochastic flow of distinguished bridges, as defined below, with $M$-coalescents.
\begin{Def} \label{flow} A flow of distinguished bridges is a collection $(B_{s,t}, -\infty<s\leq t<\infty)$ of distinguished bridges such that
\begin{itemize} 
\item $\forall s<t<u, B_{s,u}=B_{s,t}\circ B_{t,u}$ a.s. 
\item The law of $B_{s,t}$ depends only on $t-s$, and for any $s_{1}<...<s_{n}$, $B_{s_{1},s_{2}},...,B_{s_{n-1},s_{n}}$ are independent.  
\item $B_{0,0}=$Id, $B_{0,t} \rightarrow$ Id in probability when $t \rightarrow 0$.
\end{itemize}
\end{Def}
\subsubsection{Poissonian construction of distinguished flows encoding $M$-coalescent}
Let $M^{0}$ and $M^{1}$ be two independent Poissonian measures on $\mathbb{R}\times[0,1]$ with intensities $dt\otimes\nu_{0}(dx)$ and $dt\otimes\nu_{1}(dx)$. We suppose $\nu_{0}([0,1])+\nu_{1}([0,1])< \infty$ so that $N_{s,t}:=(M^{0}+M^{1})(]s,t]\times[0,1])$ is finite and $(N_{t})=(N_{0,t})_{t\geq 0}$ is a Poisson process. Let $(t^{0}_{i},x^{0}_{i})$, $(t^{1}_{i},x^{1}_{i})$  be the atoms of $M^{0}$ and $M^{1}$ in $]s,t]\times[0,1]$, we define $B_{s,t}=b_{x_{1}}\circ...\circ b_{x_{K}}$ where $K=N_{s,t}$ and where $b_{x}$ denotes $b_{0,x}$ or $b_{x,0}$ depending on whether $x$ is an atom of $M^{0}$ or $M^{1}$. From the independence of $M^{0}$ and $M^{1}$ and the independence of $M^{i}(A)$ and $M^{i}(B)$, $i=0,1$ for $A$ and $B$ disjoint, $(B_{s,t})_{s\leq t}$ is a flow in the sense of definition \ref{flow}. 
\begin{prop} 
The process $(\Pi^{0}(t), t\geq 0 )$ : $i \sim j \Leftrightarrow B^{-1}_{0,t}(U_{i})=B^{-1}_{0,t}(U_{j})$ where $B_{0,t}=b_{x_{1}}\circ...\circ b_{x_{N_{t}}}$ is a $M$-coalescent with $M=(x\nu_{0}(dx), x^{2}\nu_{1}(dx))$.
\end{prop}
\textit{Proof.} Lemma \ref{key} implies that the process $(\Pi^{0}(t))_{t\geq 0}$ corresponds to that built explicitly in Proposition \ref{construction}. $\square$
\\
The next result defines stochastic flows for general measures $\nu_{0}$ and $\nu_{1}$ on $[0,1]$. 
\begin{thm}
Let $(\nu^{n}_{0})$ and $(\nu^{n}_{1})$ be two sequences of finite measures on $[0,1]$. We call $(B^{(n)}_{s,t}, -\infty<s\leq t<\infty)$ the associated flow of bridges. Assume the weak convergences of $\Lambda^{n}_{1}(dx):=x^{2}\nu^{n}_{1}(dx)$ to $\Lambda_{1}(dx):=c_{1}\delta_{0}(dx) + x^{2}\nu_{1}(dx)$ and $\Lambda^{n}_{0}(dy):=y\nu^{n}_{0}(dy)$ to $\Lambda_{0}(dy):=c_{0}\delta_{0}+ y\mu(dy)$. We get
\begin{itemize}
\item $(B^{(n)}_{s,t}, \infty<s\leq t<\infty)$ converges, in the sense of convergence of finite-dimensional distributions, to $(B_{s,t}, \infty<s\leq t<\infty)$ a stochastic flow.
\item The process $(\Pi^{0}(t), t\geq 0)$ defined by $\Pi^{0}(t): i\sim j \Leftrightarrow B_{0,t}^{-1}(U_{i})= B_{0,t}^{-1}(U_{j})$ is a $M$-coalescent with rates $(\Lambda_{0}, \Lambda_{1})$.
\end{itemize}
\end{thm}
\textit{Proof.} We denote by $\Pi^{(n)}(s,t)$ the random partition encoded by $B^{(n)}_{s,t}$. Under the previous assumptions on $\nu_{0}$ and $\nu_{1}$, the jump rates $\lambda^{(n)}_{b,k}=\int_{0}^{1}x^{k}(1-x)^{b-k}\nu_{1}^{n}(dx)$ and $r^{(n)}_{b,k}=\int_{0}^{1}y^{k}(1-y)^{b-k}\nu^{n}_{0}(dy)$ converge: 
$$\lambda^{(n)}_{b,k} \underset{n \rightarrow \infty}{\rightarrow} c_{1}1_{k=2} + \int_{0}^{1}x^{k}(1-x)^{b-k}\nu_{1}(dx):=\lambda_{b,k}.$$
$$r^{(n)}_{b,k} \underset{n \rightarrow \infty}{\rightarrow} c_{0}1_{k=1} + \int_{0}^{1}x^{k}(1-x)^{b-k}\nu_{0}(dx):=r_{b,k}.$$ The sequence of Markov chains $(\Pi^{(n)}_{|\cro{k}}(s,t))_{t\geq s}$ converges in the sense of finite-dimensional distributions to a distinguished coalescent chain, say $\Pi^{0}_{|\cro{k}}(s,t)_{t\geq s}$. By compatibility, this implies the convergence of finite-dimensional distributions of $(\Pi^{(n)}(s,t))_{n\geq 1}$ to $\Pi^{0}(s,t)$. According to Proposition 2.9 and Lemma 4.7 in \cite{coursbertoin} (which are easily adapted to our setting), we obtain the convergence of the distinguished mass-partitions $|\Pi^{(n)}(s,t)|^{\downarrow} \underset{n \rightarrow \infty}{\rightarrow} |\Pi^{0}(s,t)|^{\downarrow}$ and the convergence of the bridge $B^{(n)}_{s,t}$ (which has jumps of size $|\Pi^{(n)}(s,t)|^{\downarrow}$) to a bridge $B_{s,t}$ (which has jumps of size $|\Pi^{0}(s,t)|^{\downarrow}$) for all $s, t \geq 0$ fixed. Thanks to the independence of $B^{(n)}_{s_{1},s_{2}},...,B^{(n)}_{s_{k-1},s_{k}}$ for any $s_{1}<...<s_{k}$ and the flow property $(B^{(n)}_{s,t}\circ B^{(n)}_{t,u}=B^{(n)}_{s,u})$, the one-dimensional convergence in distribution readily extends to finite-dimensional distribution. The existence of the flow $B$ is ensured by Kolmogorov's Extension Theorem.$\square$
\begin{rem} As mentioned in Section 2, we could define coalescents with several distinguished blocks. In particular, considering distinguished bridges which jump at $0$ and $1$, we get a flow coding a coalescent with two distinguished blocks and a population with two immigration sources.
\end{rem}
The composition of two distinguished bridges may be interpreted as the succession of two events (reproduction or immigration) in the population. A duality method provides a continuous population model.
\subsection{The dual distinguished flow and a population model with immigration}
In the same spirit of \cite{coursbertoin} and \cite{LGB1}, we interpret the dual flow, $(\hat{B}_{s,t}):=(B_{-t,-s})$ in terms of a natural model population on $[0,1[$ with fixed size $1$. We denote by $Z^{0}_{t}(dr)$, the random Stieljes measure of $\hat{B}_{0,t}$:  $Z^{0}_{t}=d\hat{B}_{0,t}$, it defines a Markov process with values in the space of probability measures on $[0,1[$ (denoted by $\mathcal{M}_{1}$). We may think of $Z^{0}_{t}(dr)$ and $Z^{0}_{t}(\{0\})$ respectively as the size of the progeny at time $t$ of the fraction $dr$ of the initial population, and as the size at time $t$ of the immigrants descendants. 
\\
The cocycle identity $\hat{B}_{t,u}\circ\hat{B}_{s,t}=\hat{B}_{s,u}$ ensures that $(Z^{0}_{t})$ is a continuous-time Markov chain with the following dynamics, whenever the measures $\nu_{0}$ and $\nu_{1}$ are finite: if $t$ is a jump time for $Z^{0}_{.}$, then the conditional law of $Z^{0}_{t}$ given $Z^{0}_{t-}$ is that of 
\begin{itemize} 
\item $(1-X)Z^{0}_{t-}+X\delta_{U}$, if $t$ is an atom of $M^{1}$, where $X$ is distributed as $\nu_{1}(.)/\nu_{1}([0,1])$ and $U$ as $Z^{0}_{t-}$.
\item $(1-Y)Z^{0}_{t-}+Y\delta_{0}$, if $t$ is an atom of $M^{0}$, where $Y$ is distributed as $\nu_{0}(.)/\nu_{0}([0,1])$.
\end{itemize}
At a reproduction time (meaning an atom of $M^{1}$) an individual picked at random in the population at generation $t-$ generates a proportion $X$ of the population at time $t$, as for the genuine generalized Fleming-Viot. At an immigration time (meaning an atom of $M^{0}$) the individual $0$ at the time $t-$ generates a proportion $Y$ of the population at time $t$. In both cases, the rest of the population at time $t-$ is reduced by a factor $1-X$ or $1-Y$ so that, at time $t$, the total size is still $1$. We call this measure-valued process a generalized Fleming-Viot process with immigration (GFVI). The genealogy of this population (which is identified as $[0,1]$) coincides with an $M$-coalescent. 
Plainly, the generator of $(Z^{0}_{t},t\geq 0)$ is $$\mathcal{L}G(\rho)=\int \nu_{1}(dx)\int \rho(da)[G((1-x)\rho+x\delta_{a})-G(\rho)] + \int \nu_{0}(dy)[G((1-y)\rho+y\delta_{0})-G(\rho)].$$
Thus, for any bounded function $G$ on $\mathcal{M}_{1}$, the space of probability measures on $[0,1]$
\begin{eqnarray*}
G(Z^{0}_{t})&-&\int_{0}^{t}ds\int \nu_{1}(dx)\int Z^{0}_{s}(da)[G((1-x)Z^{0}_{s}+x\delta_{a})-G(Z^{0}_{s})]\\
&& \qquad -\int_{0}^{t}ds\int \nu_{0}(dy)[G((1-y)Z^{0}_{s}+y\delta_{0})-G(Z^{0}_{s})]
\end{eqnarray*}
is a martingale. Considering the functions of the form $$G_{f}: \rho \in \mathcal{M}_{1} \mapsto \int_{[0,1]^{p}}f(x_{1},...,x_{p})\rho(dx_{1})...\rho(dx_{p})=\langle f,\rho^{\otimes p} \rangle,$$ for $f$ a continuous function on $[0,1]^{p}$, we generalize in the following lemma this result for infinite measures.
\begin{lem} \label{martingale} Assume that $\nu_{0}$ and $\nu_{1}$ have infinite masses and $c_{0}$, $c_{1}$ are zero, we define the operator $\mathcal{L}$, acting on functions of the type $G_{f}$, by
$$\mathcal{L}G_{f}(\rho)=\int \nu_{1}(dx)\int \rho(da)[G_{f}((1-x)\rho+x\delta_{a})-G_{f}(\rho)] + \int \nu_{0}(dy)[G_{f}((1-y)\rho+y\delta_{0})-G_{f}(\rho)].$$
The process $G_{f}(Z^{0}_{t})-\int_{0}^{t}\mathcal{L}G_{f}(Z^{0}_{s})ds$ is a martingale.
\end{lem}
\textit{Proof.} Consider two sequences $(\nu_{1}^{n})$ and $(\nu_{0}^{n})$ of finite measures on $[0,1[$, suppose that $\Lambda^{n}_{1}(dx):= x^{2}\nu_{1}^{n}(dx)$ and $\Lambda^{n}_{0}(dx):=x\nu^{n}_{0}(dx)$ weakly converge to some finite measures $\Lambda_{1}(dx)$ and $\Lambda_{0}(dx)$.
\\
For $G_{f}(\rho)=\prod_{i}^{p}<\rho,\phi_{i}>$, Bertoin and Le Gall obtain in \cite{LGB1}: 
$$\int \nu_{1}^{n}(dx)\int \rho(da)[G((1-x)\rho+x\delta_{a})-G(\rho)]=\sum_{\substack{I \subset \{1,...,p\}\\ |I|\geq 2}}\lambda^{n}_{p,|I|}\int \rho(dx_{1})...\rho(dx_{p})[f(x^{I}_{1},...,x^{I}_{p})-f(x_{1},...,x_{p})].$$ 
With $(x^{I}_{1},...,x^{I}_{p})=(y_{1},...,y_{p})$ where for all $i \in I$, $y_{i}=x_{\inf I}$ and the values $y_{i}$, $i\notin I$ listed in the order of $\{1,...,p\}\setminus I$ are the numbers $x_{1},..,x_{\inf I-1},x_{\inf I+1},...,x_{p-|I|+1}$.\\
\\
The assumption on $\nu_{1}$ ensures that the right side converges to $$\sum_{I \subset \{1,...,p\}; |I|\geq 2}\lambda_{p,|I|}\int \rho(dx_{1})...\rho(dx_{p})[f(x^{I}_{1},...,x^{I}_{p})-f(x_{1},...,x_{p})].$$
It remains to study the "immigration" part, that is to establish
$$\int \nu^{n}_{0}(dy) [G((1-y)\rho+y\delta_{0})-G(\rho)] =\sum_{J \subset \{1,...,p\}; |J|\geq 1} r^{n}_{p,|J|}\int \rho(dx_{1})...\rho(dx_{p})[f(x^{0}_{1},...,x^{0}_{p})-f(x_{1},...,x_{p})]$$
with $(x^{0}_{1},...,x^{0}_{p})=(z_{1},...,z_{p})$ where for all $i \in J$, $z_{i}=0$ and the values $z_{i}$, $i\notin J$ listed in the order of $\{1,...,p\}\setminus J$ are the numbers $x_{1},..,x_{\inf J-1},x_{\inf J+1},...,x_{p-|J|}$.
\\
\\
An easy calculation gives
$$G((1-y)\rho+y\delta_{0})=\prod_{i=1}^{p}[(1-y)\langle \rho,\phi_{i}\rangle+y\phi_{i}(0)]=\sum_{J\subset \{1,...,p\}}(1-y)^{p-|J|}y^{|J|}\prod_{j\notin J}\langle\rho,\phi_{j}\rangle \prod_{j\in J}\phi_{j}(0)$$
and then from the following identity $$\prod_{j\notin J}\langle\rho,\phi_{j}\rangle\prod_{j\in J}\phi_{j}(0)=\int_{[0,1]^{p}}f(x^{0}_{1},...,x^{0}_{p})\rho(dx_{1})...\rho(dx_{p})$$ it follows that
$$\int \nu^{n}_{0}(dy) G((1-y)\rho+y\delta_{0})=\sum_{J \subset \{1,...,p\}; |J|\geq 1} r^{n}_{p,|J|}\int \rho(dx_{1})...\rho(dx_{p})f(x^{0}_{1},...,x^{0}_{p}).$$
Moreover, the right hand side converges to $$\sum_{J \subset \{1,...,p\}; |J|\geq 1} r_{p,|J|}\int \rho(dx_{1})...\rho(dx_{p})f(x^{0}_{1},...,x^{0}_{p}),$$ by passing to the limit in $n\rightarrow \infty$, it follows that for $f(x_{1},...,x_{p})=\prod_{i}^{p}\phi_{i}(x_{i})$, the process $$M_{f}(t):=G_{f}(Z^{0}_{t})-\int_{0}^{t}LG_{f}(Z^{0}_{s})ds$$ is a martingale where $L$ is the operator defined by
\begin{align*}
LG_{f}(\rho)&=\sum_{I \subset \{1,...,p\}; |I|\geq 2} \lambda_{p,|I|} \! \int \rho(dx_{1})...\rho(dx_{p})[f(x^{I}_{1},...,x^{I}_{p})-f(x_{1},...,x_{p})]\\
&+ \sum_{J\subset\{1,...,p\}; |J|\geq 1} r_{p,|J|} \! \int \rho(dx_{1})...\rho(dx_{p})[f(x^{0}_{1},...,x^{0}_{p})-f(x_{1},...,x_{p})].
\end{align*} 
Since any continuous function on $[0,1]^{p}$ is the uniform limit of linear combinations of functions of the previous type, we easily conclude that $M_{f}$ is a martingale for any continuous function on $[0,1]^{p}$. The statement claims that when $c_{0}=c_{1}=0$ the generator has an integral form as the one obtained for finite measures. We assume now that $c_{0}$ and $c_{1}$ are zero. Let $A_{1},...,A_{p}$ i.i.d variables distributed as $\rho$, and $x,y\in [0,1]$. Let $(\beta_{j})$, $(\beta'_{j})$ be two sequences of Bernoulli variables of parameters $x$ and $y$. We set $I:=\{j, \beta_{j}=1\}$ and $J:=\{j, \beta'_{j}=1\}$. Let $f$ be a continuous function on $[0,1]^{p}$, for $G_{f}(\rho)=\langle\rho^{\otimes p},f\rangle$, it is readily checked (see \cite{Birk}) that 
\begin{align*}
&\int \rho(da)[G_{f}((1-x)\rho +x\delta_{a})-G_{f}(\rho)]= \mathbb{E}[f(A_{1}^{J}...,A_{p}^{J})]-\mathbb{E}[f(A_{1},...,A_{p})]\\
&=\sum_{I\subset \{1,...,p\}, |I|\geq 2} \! x^{|I|}(1-x)^{p-|I|}\int \rho(dx_{1})...\rho(dx_{p})(f(x^{I}_{1},...,x^{I}_{p})-f(x_{1},...,x_{p})),\\
\\
&G_{f}((1-y)\rho+y\delta_{0})-G_{f}(\rho)= \mathbb{E}[f(A_{1}^{0}...,A_{p}^{0})]-\mathbb{E}[f(A_{1},...,A_{p})]\\
&=\sum_{J\subset \{1,...,p\}, |J|\geq 1} \! y^{|J|}(1-y)^{p-|J|}\int \rho(dx_{1})...\rho(dx_{p})(f(x^{0}_{1},...,x^{0}_{p})-f(x_{1},...,x_{p})).
\end{align*}
We deduce that the process $(Z^{0}_{t})$ solves the following martingale problem: for any continuous function $f$ on $[0,1]^{p}$, $G_{f}(Z^{0}_{t})-\int_{0}^{t}ds\mathcal{L}G_{f}(Z^{0}_{s})$ is a martingale. $\square$
\begin{prop} 
The law of the process $(Z^{0}_{t}, t\geq 0)$ is characterized by the martingale problem of Lemma \ref{martingale}, and the operator $\mathcal{L}$ is an extended generator of the process $(Z^{0}_{t}, t\geq 0)$.
\end{prop}
\textit{Proof.} We will use the same duality argument as in Bertoin and Le Gall \cite{LGB1}. With their notation, we define a class of functions from $\mathcal{M}_{1}\times \mathcal{P}^{0}_{p}$ to $\mathbb{R}$ 
\begin{center}
$\Phi_{f} : (m,  \pi) \in \mathcal{M}_{1}\times \mathcal{P}^{0}_{p} \mapsto \int \delta_{0}(dx_{0})m(dx_{1})...m(dx_{\#\pi-1})f(Y(\pi; x_{1},...,x_{\#\pi-1}))$
\end{center}
with $f$ a continuous function on $[0,1]^{p}$ and $Y(\pi; x_{1},...,x_{\#\pi-1})=(y_{1},...,y_{p})$ such that $y_{j}=x_{i}$ if $ j \in \pi_{i}$ for any $i\geq 0$.
\\
For a fixed partition $\pi$ in $\mathcal{P}^{0}_{n}$, there exists a function $g$ continuous on $[0,1]^{\#\pi-1}$ with $\mu \mapsto \Phi_{f}(\mu,\pi)=G_{g}(\mu)$ and then $\mathcal{L}\Phi_{f}(\mu,\pi)=\mathcal{L}G_{g}(\mu)$ is well-defined. We stress that for a fixed measure $\mu$, $\Phi_{f}(\mu,.)$ is a function on $\mathcal{P}^{0}_{p}$.
We show the following duality result: $$\mathbb{E}[\Phi_{f}(Z^{0}_{0},\Pi^{0}(t))]=\mathbb{E}[\Phi_{f}(Z^{0}_{t},\Pi^{0}_{0})].$$ 
By the cocycle property of the stochastic flow involved, it suffices to focus on process beginning at $Z^{0}_{0}=\lambda$.
\begin{align*}
\mathbb{E}_{\lambda}[\phi_{f}(Z^{0}_{t},0_{[p]})]&=\mathbb{E}[\int_{[0,1]^{p+1}}\delta_{0}(dx_{0})d\hat{B}_{0,t}(x_{1})...d\hat{B}_{0,t}(x_{p})f(x_{1},...,x_{p})]\\
&=\mathbb{E}[\int_{[0,1]^{p+1}}\delta_{0}(dx_{0})dx_{1}...dx_{p}f(\hat{B}^{-1}_{0,t}(x_{1}),...,\hat{B}^{-1}_{0,t}(x_{p}))]\\
&=\mathbb{E}[f(\hat{B}^{-1}_{0,t}(V_{1}),...,\hat{B}^{-1}_{0,t}(V_{p}))]
\end{align*}
where $(V_{i}, 1\leq i\leq p)$ are independent and uniformly distributed on $[0,1]$. We define for $1\leq i\leq \#\Pi^{0}_{|\cro{p}}(t)-1$, $V'_{i}:=\hat{B}^{-1}_{0,t}(V_{j})$ for $j\in \Pi^{0}_{i|\cro{p}}(t)$. By Lemma \ref{key}, $(V'_{i})_{1\leq i\leq\#\Pi^{0}_{|\cro{p}}(t)-1}$ are uniform iid, independent of $\Pi^{0}(t)$ where $\Pi^{0}$ is a $M$-coalescent with rates $(x\nu_{0}(dx),x^{2}\nu_{1}(dx))$.
\\
We get, $$\mathbb{E}_{\lambda}[\phi_{f}(Z^{0}_{t},0_{\cro{p}})]=\mathbb{E}[\int_{[0,1]^{p+1}}\delta_{0}(dx_{0})dx_{1}...dx_{\#\Pi^{0}_{|\cro{p}}(t)-1}
f(y_{1},...y_{\#\Pi_{|\cro{p}}^{0}(t)-1})]$$ with $y_{j}=x_{i}$ if $j \in \Pi^{0}_{i|\cro{p}}(t)$.
\\ 

Thus, we deduce the duality result:
\begin{center}
$\mathbb{E}_{\lambda}[\Phi_{f}(Z^{0}_{t},0_{\cro{p}})]=\mathbb{E}[\Phi_{f}(\lambda,\Pi^{0}_{|\cro{p}}(t))]$
and then
$\mathcal{L}\Phi_{f}(\mu,\pi)=\mathcal{L}^{*}\Phi_{f}(\mu,\pi).$
\end{center}
From Theorem 4.4.2 in \cite{EthierKurtz}, this implies uniqueness for the martingale problem, as well as strong Markov property for the solution.$\square$ 
\begin{rem} In the case of a standard $M$-coalescent, $Z^{0}_{0}$ is the Lebesgue measure $\lambda$, and we have
\begin{center}
$Z^{0}_{t}(dr)=|\Pi^{0}_{0}(t)|\delta_{0}(dr)+\sum_{i\geq 1}|\Pi^{0}(t)|^{\downarrow}_{i}\delta_{W_{i}}(dr)+(1-\sum_{i\geq 0}|\Pi^{0}(t)|^{\downarrow}_{i})dr$
\end{center}
where $(W_{i}, i\geq 1)$ are independent uniform and independent of $\Pi^{0}(t)$.  
\end{rem}
The extinction of the initial types corresponds to the absorption of the GFVI process $(Z^{0}_{t}, t\geq 0)$ at $\delta_{0}$. Plainly this event occurs if and only if the measure $\Lambda_{0}$ is not the zero measure and the $M$-coalescent embedded is coming down from infinity. By Theorem \ref{cdi}, we know that the coming down from infinity depends only on the measure $\Lambda_{1}$. In terms of the population model, the immigration mechanism, encoded by $\Lambda_{0}:=c_{0}\delta_{0}+x\nu_{0}(dx)$, has no impact on the extinction occurrence, provided of course that $\Lambda_{0}$ is not the zero measure.
\section{Proof of Theorem \ref{cdi}}
We recall the statement of Theorem \ref{cdi} and give a proof based on martingale arguments.\\ \\
\textbf{ Theorem \ref{cdi} } \textit{ The $(\Lambda_{0},\Lambda_{1})$-coalescent comes down from infinity if and only if 
\begin{center} $(\Lambda_{0}+\Lambda_{1})(\{1\})>0$ or $\sum_{n=2}^{\infty}\frac{1}{\phi_{1}(n)}<\infty$ \end{center} where $\phi_{1}(n)=\sum_{k=2}^{n}(k-1)C^{k}_{n}\lambda_{n,k}$ and $\lambda_{n,k}=\int_{0}^{1}x^{k-2}(1-x)^{n-k}\Lambda_{1}(dx)$.}
\\
\\
For the $\Lambda$-coalescents, (in our setting it corresponds to have $\Lambda_{0}\equiv 0$), Schweinsberg studied the mean time of "coming down from infinity" and conclude using the Kochen-Stone lemma. The proof we give here is based on martingale arguments. To show that the convergence of the series is sufficient for the coming down from infinity, we need to prove Lemma \ref{CS} which is close to the one of Proposition 4.9 page 202 in \cite{coursbertoin}. The necessary part of the proof does not follow Schweinsberg's ideas. Assuming that the coalescent comes down from infinity and the sum is infinite, we will define a supermartingale (thanks to Lemmas \ref{binom}, \ref{half} and \ref{superm}) and find a contradiction (Lemma \ref{CN}).
\begin{lem} \label{CS}
Let $(\Pi^{0}(t),t\geq 0)$ be a $(\Lambda_{0}, \Lambda_{1})$-coalescent where $\Lambda_{0}+\Lambda_{1}$ has no mass at $1$. Let us define the fixation time $$\zeta:=\inf\{t\geq 0, \Pi^{0}(t)=\{\mathbb{Z}_{+},\emptyset,...\}\}.$$ 
Define $$\phi(n)=\sum_{k=2}^{n}(k-1)C^{k}_{n}\lambda_{n,k}+\Lambda_{0}([0,1])n,$$ 
then the expectation of fixation time is bounded by $$\mathbb{E}[\zeta]\leq \sum_{n=1}^{\infty}1/\phi(n).$$
As a consequence, if the series in the right-hand side converges, the fixation time is finite with probability one.
\end{lem}
\textit{Proof.} We shall study the process of blocks which do not contain $0$: for all $t >0$, $\Pi^{*}(t):=\{\Pi^{0}_{1}(t),...\}$. This process is not partition-valued. The jump rates of $\#\Pi^{*}_{|\cro{n}}(t)$ are easily computed: for $2 \leq k\leq l+1$, $\#\Pi^{*}_{|\cro{n}}$ jumps from $l$ to $l-k+1$ with rate: $$C^{k}_{l}\lambda_{l,k}1_{k\leq l}+C^{k-1}_{l}r_{l,k-1}.$$ The first term represents coagulation of $k$ blocks not containing $0$ in one, the second represents disappearance of $k-1$ blocks (coagulation with $\Pi^{0}_{0}$). We get the infinitesimal generator of $\#\Pi^{*}_{|\cro{n}}$ :
$$G^{[n]}f(l)=\sum_{k=2}^{l}C^{k}_{l}\lambda_{l,k}[f(l-k+1)-f(l)]+\sum_{k=2}^{l+1}C^{k-1}_{l}r_{l,k-1}[f(l-k+1)-f(l)].$$
We define \begin{center} $\phi_{1}(n)=\sum_{k=2}^{n}(k-1)C^{k}_{n}\lambda_{n,k}$ and $\phi_{2}(n)=\sum_{k=2}^{n+1}(k-1)C^{k-1}_{n}r_{n,k-1}$. \end{center} Using the binomial formula, we get $\phi_{2}(n)=\Lambda_{0}([0,1])n$.  We remark that $\phi_{1}$ is an increasing function. Setting $\phi(n)=\phi_{1}(n)+\phi_{2}(n)$ and assuming the convergence of the sum $\sum_{n=1}^{\infty}1/\phi(n)$, we define $$f(l)=\sum_{k=l+1}^{\infty}1/\phi(k).$$ The map $\phi$ is increasing, we thus have $f(l-k+1)-f(l)\geq \frac{k-1}{\phi(l)}$ and then $$G^{[n]}f(l)\geq \sum_{k=2}^{l+1}(C^{k}_{l}\lambda_{l,k}+C^{k-1}_{l}r_{l,k-1})\frac{k-1}{\phi(l)}=1$$ 
The process $f(\#\Pi^{*}_{|\cro{n}}(t))-\int_{0}^{t}G^{[n]}f(\#\Pi^{*}_{|\cro{n}}(s))ds$ is a martingale. The quantity $$\zeta_{n}:=\inf\{t; \#\Pi^{*}_{\cro{n}}(t)=0\}$$ is a finite stopping time. Let $k\geq 1$, applying the optional stopping theorem to the bounded stopping time $\zeta_{n}\wedge k$, we get : $$\mathbb{E}[f(\#\Pi^{*}_{|\cro{n}}(\zeta_{n}\wedge k))]-\mathbb{E}[\int_{0}^{\zeta_{n}\wedge k}G^{[n]} f(\#\Pi^{*}_{|\cro{n}}(s))ds]=f(n)$$
With the inequality $G^{[n]}f(l)\geq 1$, we deduce that $$\mathbb{E}[\zeta_{n}\wedge k]\leq \mathbb{E}[f(\#\Pi^{*}_{|\cro{n}}(\zeta_{n}\wedge k))]-f(n).$$ By monotone convergence and Lebesgue's theorem, we have $\mathbb{E}[\zeta_{n}]\leq f(0)-f(n)$.  Passing to the limit in $n$, we have $\zeta_{n} \uparrow \zeta_{\infty}:= \inf \{t; \#\Pi(t)=1\}$ and $f(n) \longrightarrow 0$, thus $$\mathbb{E}[\zeta_{\infty}]\leq f(0)=\sum_{k=1}^{\infty}1/\phi(k).$$ 
$\square$
\\
\\
By simple series comparisons, we deduce the sufficient part of the theorem. Plainly $\phi(n) \geq \phi_{1}(n)$, and if the series $\sum_{n=1}^{\infty}1/\phi_{1}(n)$ converges, then by Lemma \ref{CS} the $M$-coalescent comes down from infinity.
\\
\\
To show that the convergence of the series is necessary for the coming down, we must look more precisely at the behavior of jumps. The next technical lemmas show that when a distinguished coalescent comes down from infinity, there is a finite number of jumps which make decrease by half or more the number of blocks. Lemma \ref{half} will allow us to study the process of blocks number before the first of these times. Assuming that the sum $\sum_{n=2}^{\infty}\frac{1}{\phi(n)}=\infty$ is infinite, we will define a supermartingale in Lemma 17 and find a contradiction by applying the optional stopping theorem. 
\\
\\
We have already seen that, as for the $\Lambda$-coalescent, a way to understand the dynamics of a $(\Lambda_{0}, \Lambda_{1})$-coalescent, when $\Lambda_{0}, \Lambda_{1}$ have no mass at $0$, is to imagine drawing an infinite sequence of Bernoulli variables at each jump time, with parameter $x$ controlled by the measures $\nu_{0}(dx)=x^{-1}\Lambda_{0}(dx)$ and $\nu_{1}(dx)=x^{-2}\Lambda_{1}(dx)$. The following technical lemma allows to estimate the chance for a Bernoulli vector to have more than half terms equals to $1$. 
\begin{lem} \label{binom}
Let $(X_{1},X_{2},...)$ be independent Bernoulli variables with parameter $x\in [0,1/4[$. Defining $S^{(x)}_{n}=X_{1}+...+X_{n}$, for every $n_{0}$, there is the bound
$$\mathbb{P}[\exists n\geq n_{0}; S^{(x)}_{n}>\frac{n}{2}]\leq \frac{\exp{(-n_{0}f(x))}}{1-\exp{(-f(x))}}$$
with $f(x) \sim \frac{1}{2}log(1/x)$ when $x \rightarrow 0$.
\end{lem}
\textit{Proof}. By Markov inequality, for all $t>0$, $$\mathbb{P}[S_{n}^{(x)}\geq n/2]\leq e^{-nt/2} \mathbb{E}[e^{tS^{(x)}_{n}}]=\exp{(-n[t/2-\log(xe^{t}+1-x)])}.$$ Applying this inequality for $t=\log(1/x)$, we get $\mathbb{P}[S_{n}^{(x)}\geq n/2]\leq e^{-nf(x)}$ where $$f(x)=\frac{1}{2}\log(1/x)-\log(2-x).$$ The function $f$ is non-negative on $]0,\frac{1}{4}[$ and then we obtain the convergence of the geometric sum
$$\mathbb{P}[\exists n\geq n_{0}; S^{(x)}_{n}>\frac{n}{2}]\leq \sum_{n\geq n_{0}}\exp{(-nf(x))}=\frac{\exp{(-n_{0}f(x))}}{1-\exp{(-f(x))}}.$$
Moreover we have $f(x) \underset{x\rightarrow 0^{+}}{\sim}\frac{1}{2}\log(1/x)$.
\begin{lem} \label{half} Assume that the $M$-coalescent comes down from infinity. With probability one, we have
$$\tau:=\inf\{t>0, \#\Pi^{0}(t)<\frac{\#\Pi^{0}(t-)}{2}\}>0.$$ Moreover, if we define $\tau_{n}:=\inf\{t>0, \#\Pi^{0}_{|\cro{n}}(t)<\frac{\#\Pi^{0}_{|\cro{n}}(t-)}{2}\}$, then the sequence of stopping times $\tau_{n}$ converges to $\tau$ almost surely.
\end{lem}
\textit{Proof.} 
Obviously, binary coagulations play no role in the statement and we may assume that $\Lambda_{0}(\{0\})=0$ and $\Lambda_{1}(\{0\})=0$. Let $N$ be a Poisson measure with intensity $dt\otimes(x^{-1}\Lambda_{0}(dx)+x^{-2}\Lambda_{1}(dx))$. Recall the notation in Lemma \ref{binom}. Let $n_{0}\geq 4$, we will show that
$$N(\{(t,x); t\leq 1; \exists n\geq n_{0}; S^{(x)}_{n}\geq n/2\})<\infty.$$
By Proposition \ref{coming}, $\#\Pi^{0}(\epsilon)<\infty$ a.s for every $\epsilon>0$. We will then deduce that there is a finite number of jump times before $1$ where more than half blocks coagulate. By the Feller property and therefore the regularity of paths of $(\Pi^{0}(t))_{t\geq 0}$, $\Pi^{0}(0+)=\Pi^{0}(0)$ and $0$ is not a jump time, then almost surely: $\tau>0$. Moreover, $\#\Pi^{0}(\tau-)<\infty$ and then for all $n\geq \#\Pi^{0}(\tau-)$, $\tau_{n}=\tau$. We deduce that $\tau_{n}\underset{n\rightarrow \infty}{\rightarrow} \tau$ almost surely.
\\
By Poissonian calculations, we get $$\mathbb{E}[N(\{(t,x); t\leq 1; \exists n\geq n_{0}, S^{(x)}_{n}>n/2\})]=\int_{0}^{1}(\nu_{0}+\nu_{1})(dx)\mathbb{P}[\exists n\geq n_{0}, S^{(x)}_{n}>n/2].$$
By Lemma \ref{binom}, we get
$$\int_{0}^{1}(\nu_{0}+\nu_{1})(dx)\mathbb{P}[\exists n\geq n_{0}, S^{(x)}_{n}>n/2]\leq \int_{0}^{\frac{1}{4}}(\nu_{0}+\nu_{1})(dx)\frac{\exp{(-n_{0}f(x))}}{1-\exp{(-f(x))}}+\int_{1/4}^{1}(\nu_{0}+\nu_{1})(dx).$$
On the one hand $$\int_{0}^{\frac{1}{4}}(\nu_{0}+\nu_{1})(dx)\frac{\exp{(-n_{0}f(x))}}{1-\exp{(-f(x))}}<\infty$$ because the integrand is bounded by $8x^{\frac{n_{0}}{2}}$ and $n_{0}\geq 4$,
\\  
\\
on the other hand $$\int_{1/4}^{1}(\nu_{0}+\nu_{1})(dx)\leq \int_{1/4}^{1}x^{-2}x^{2}(\nu_{0}+\nu_{1})(dx)\leq 16\int_{0}^{1}x^{2}(\nu_{0}+\nu_{1})(dx)<\infty.$$ This completes the proof.$\square$
\\
\\
Assuming that the coalescent comes down from infinity and that $\sum_{n\geq 1} \frac{1}{\phi(n)} =\infty$, we can define a supermartingale. We will find a contradiction using the optional stopping theorem.\\ 
We define the decreasing function: $$f(n)=\exp(-\sum_{k=1}^{n+1}\frac{1}{\phi(k)})$$ where $\phi(n)=\phi_{1}(n)+\phi_{2}(n)$ with $\phi_{1}(n)$ and $\phi_{2}(n)$ are defined as in Lemma \ref{CS}. 
\begin{lem} \label{superm}
There exists a constant $C>0$ such that for all $n\geq 1$, $(e^{-Ct}f(\#\Pi^{*}_{|[n]}(t)))_{t\leq \tau_{n}}$ is a non-negative supermartingale.
\end{lem}
\textit{Proof.} We recall that the generator of $(\#\Pi^{*}_{|\cro{n}}(t))_{t\geq 0}$ is $$G^{[n]}g(l)=\sum_{k=2}^{l+1} [C^{k}_{l}1_{k\leq l}\lambda_{l,k}+C^{k-1}_{l}r_{l,k-1}][g(l-k+1)-g(l)].$$ Stopping the process at $\tau_{n}$, the jump times where more half of blocks coagulate are ignored, and the generator of the stopped process is
$$A^{[n]}g(l)=\sum_{k=2}^{l/2+1} [C^{k}_{l}1_{k\leq l/2}\lambda_{l,k}+C^{k-1}_{l}r_{l,k-1}][g(l-k+1)-g(l)].$$
We set $\Psi(n)=\int_{0}^{1}(e^{-nx}-1+nx)\nu_{1}(dx)$, an easy verification allows to claim the existence of $c>0$ such that $c\Psi(q)\leq \phi_{1}(q)\leq \Psi(q)$ (see remark p170 in\cite{LGB3}), moreover
$$\frac{\Psi(q)}{q}=\int_{0}^{1}(1-e^{-qx})\nu_{1}([x,1[)dx \underset{q \rightarrow \infty}\rightarrow \int_{0}^{1}x\nu_{1}(dx)>0.$$
Plainly, $h(q)=\Psi(q)/q$ is a concave function, so $h(q/2)\geq h(q)/2$ and $\Psi(q/2)\geq \Psi(q)/4$. \\
Let us compute
$$A^{[n]}f(l)=\sum_{k=2}^{l/2+1} [C^{k}_{l}1_{k\leq l/2}\lambda_{l,k}+C^{k-1}_{l}r_{l,k-1}]f(l)[\exp(\sum_{l-k+2}^{l}\frac{1}{\phi(j)})-1].$$
We have $$\sum_{l-k+2}^{l}\frac{1}{\phi(j)}\leq \frac{k-1}{\phi(l-k+2)} \leq \frac{k-1}{\phi(l/2)} \quad \forall k \leq \frac{l}{2}+1$$ and $e^{x}-1\leq cx$ for small $x$, then for large $l$
$$A^{[n]}f(l)\leq cf(l)\sum_{k=1}^{l/2+1}[C^{k}_{l}\lambda_{l,k}+C^{k-1}_{l}r_{l,k-1}]\frac{k-1}{\phi(l/2)}.$$
Thus, $A^{[n]}f(l)\leq cf(l)\frac{\phi(l)}{\phi(l/2)}$. 
\\
By definition $\phi(l)=\phi_{1}(l)+\Lambda_{0}([0,1])l$. Moreover $l/\phi_{1}(l)$ is bounded (it converges to $(\int_{[0,1]}x\nu_{1}(dx))^{-1}$) and from the inequalities: $c\Psi(l)/l\leq \phi_{1}(l)/l\leq \Psi(l)/l$ and $\Psi(l/2)\geq \Psi(l)/4$, we deduce that for some constant $C>0$ 
$$\frac{\phi(l)}{\phi(l/2)}\leq \frac{\phi_{1}(l)+\Lambda([0,1])l}{\phi_{1}(l/2)}=\frac{\phi_{1}(l)}{\phi_{1}(l/2)}+(\Lambda_{0}([0,1]))\frac{l}{\phi_{1}(l)}\frac{\phi_{1}(l)}{\phi_{1}(l/2)}\leq \frac{4}{c}[1+\Lambda_{0}([0,1])\frac{l}{\phi_{1}(l)}]\leq C.$$ Therefore, $A^{[n]}f(l)\leq Cf(l)$ and $(e^{-Ct}f(\#\Pi^{*}_{|[n]}(t)))_{t\leq \tau}$ is a supermartingale. $\square$
\begin{lem} \label{CN}
If $\sum_{n\geq 1} \frac{1}{\phi(n)} =\infty$ then $\Pi^{0}$ does not come down from infinity.
\end{lem}
\textit{Proof.} Assume that the $M$-coalescent comes down from infinity. By Proposition \ref{coming}, we know that $T=0$ a.s. Let $T^{(n)}_{j}:=\inf\{t>0; \#\Pi^{*}_{|[n]}(t)\leq j\}$. We apply to the previous supermartingale, the optional stopping theorem at time $T^{(n)}_{j}\wedge \tau_{n}$ and get
$$\mathbb{E}[e^{-cT^{(n)}_{j}\wedge\tau_{n}}f(\#\Pi^{*}_{|[n]}(\tau_{n}\wedge T^{(n)}_{j}))]\leq f(n).$$
Passing to the limit with $n \uparrow \infty$: $f(n)\rightarrow 0$, $T^{(n)}_{j} \uparrow T_{j}$ and $\tau_{n} \rightarrow \tau>0$ (by Lemma \ref{half}). The time $T_{j}$ is strictly positive for some $j$, then $\tau\wedge T_{j}>0$, $\#\Pi^{*}(\tau\wedge T_{j})<\infty$ and thus, $f(\#\Pi^{*}(\tau\wedge T_{j}))>0$ almost surely. We have 
$$\mathbb{E}[e^{-cT_{j}\wedge \tau}]=0$$
then $T_{j}=\infty$ a.s, which is not possible on $T<\infty$.$\square$
\\
\\
It remains to establish that the convergence of the series is necessary for the coming down from infinity. When $\Lambda_{0}(\{0\})=0$, the previous lemma claims that if the $M$-coalescent comes down from infinity then $\sum_{n\geq 1} \frac{1}{\phi(n)} <\infty$. It suffices to show that $$\sum_{n\geq 2}\frac{1}{\phi_{1}(n)+\Lambda_{0}([0,1])n}<\infty \Longrightarrow  \sum_{n\geq 2}\frac{1}{\phi_{1}(n)}<\infty.$$ The sequence $(\frac{\Psi(n)}{n})_{n\geq 1}$ is increasing and tends to $\int_{0}^{1}x^{-1}\Lambda_{1}(dx)$ (possibly infinite). From the inequality: $c\frac{\Psi(n)}{n}\leq \frac{\phi_{1}(n)}{n}\leq \frac{\Psi(n)}{n}$, we get that $\frac{n}{\phi_{1}(n)}$ is bounded. It follows that $$\frac{1}{\phi_{1}(n)+\Lambda_{0}([0,1])n}=\frac{1}{\phi_{1}(n)(1+\Lambda_{0}([0,1])\frac{n}{\phi_{1}(n)})}\geq c\frac{1}{\phi_{1}(n)}$$ for some constant $c>0$.\\
We then get the necessary part, and combining the results, Theorem \ref{cdi} is deduced.
\\
\\
\textbf{Acknowledgments}. This is a part of my PhD thesis. I would like to thank my advisor Jean Bertoin for introducing this subject, his useful advice, and his encouragement. I thank the referee for his very careful reading and helpful suggestions.

\end{document}